\documentclass{amsart}
\usepackage{verbatim,amsmath,amssymb,amscd,color,graphics, axodraw}
\usepackage{pstricks,pst-node}
\usepackage[all]{xy}
\usepackage[vcentermath]{myyoungtab}

\newcommand{\mb}{\text{}}
\newcommand{\ab}{\bar{1}}
\newcommand{\bb}{\bar{2}}
\newcommand{\cb}{\bar{3}}
\newcommand{\db}{\bar{4}}
\newcommand{\eb}{\bar{5}}

\newcommand{\aal}{\mathbf{\stackrel{\leftarrow}{a}}}
\newcommand{\aar}{\mathbf{\stackrel{\rightarrow}{a}}}
\newcommand{\bbr}{\mathbf{\stackrel{\rightarrow}{b}}}
\newcommand{\bij}{\Phi}
\newcommand{\cl}{\mathrm{cl}}

\newcommand{\diagram}{\mathrm{diagram}}
\newcommand{\gc}{\geh_0}
\newcommand{\geh}{\mathfrak{g}}

\newcommand{\inner}[2]{\langle #1\,,\,#2\rangle}
\newcommand{\is}{\mathrm{inner}}

\newcommand{\la}{\lambda}
\newcommand{\La}{\Lambda}
\newcommand{\lev}{\mathrm{lev}}
\newcommand{\om}{\omega}
\newcommand{\os}{\mathrm{outer}}

\newcommand{\sigD}{\mathfrak{S}}
\newcommand{\ve}{\varepsilon}
\newcommand{\vp}{\varphi}
\newcommand{\wt}{\mathrm{wt}}

\newcommand{\redplus}{\textcolor{red}{\oplus}}
\newcommand{\redminus}{\textcolor{red}{\ominus}}
\newcommand{\blueminus}{\textcolor{blue}{\ominus}}
\newcommand{\greenplus}{\textcolor{green}{\raisebox{.07cm}{\scalebox{.5}{\framebox{+}}}}}
\newcommand{\greenminus}{\textcolor{green}{\raisebox{.07cm}{\scalebox{.5}{\framebox{-}}}}}

\newcommand{\ba}[1]{\overline{#1}}

\numberwithin{equation}{section}

\newtheorem{theorem}{Theorem}
\newtheorem{proposition}[theorem]{Proposition}
\newtheorem{lemma}[theorem]{Lemma}

\theoremstyle{definition}
\newtheorem{property}{Property}
\newtheorem{definition}{Definition}
\newtheorem{remark}{Remark}
\newtheorem{example}[remark]{Example}

\numberwithin{theorem}{section}
\numberwithin{definition}{section}
\numberwithin{remark}{section}

\begin{document}

\title[Kirillov--Reshetikhin crystals of type $D_n^{(1)}, B_n^{(1)}, A_{2n-1}^{(2)}$]
{Combinatorial structure of Kirillov--Reshetikhin crystals of type $D_n^{(1)}$,
$B_n^{(1)}$, $A_{2n-1}^{(2)}$}

\author[Anne Schilling]{Anne Schilling}
\address{Department of Mathematics, University of California, One Shields
Avenue, Davis, CA 95616-8633, U.S.A.}
\email{anne@math.ucdavis.edu}
\urladdr{http://www.math.ucdavis.edu/\~{}anne}
\thanks{\textit{Date:} April 2007, revised October 2007}
\thanks{Partially supported by NSF grant DMS-0501101.}
 
\begin{abstract}
We provide the explicit combinatorial structure of the Kirillov--Reshetikhin crystals
$B^{r,s}$ of type $D_n^{(1)}$, $B_n^{(1)}$, and $A_{2n-1}^{(2)}$. This is achieved by 
constructing the crystal analogue $\sigma$ of the automorphism of the $D_n^{(1)}$ 
(resp. $B_n^{(1)}$ or $A_{2n-1}^{(2)}$) Dynkin diagram that interchanges the 0
and 1 node. The involution $\sigma$ is defined in terms of new $\pm$ diagrams that govern 
the $D_n$ to $D_{n-1}$ (resp. $B_n$ to $B_{n-1}$, or $C_n$ to $C_{n-1}$) branching.  
It is also shown that the crystal $B^{r,s}$ is perfect. These crystals have been implemented
in MuPAD-Combinat; the implementation is discussed in terms of many examples.
\end{abstract}

\maketitle

\section{Introduction}

The irreducible finite-dimensional modules over a quantized affine
algebra $U'_q(\geh)$ were classified by Chari and
Pressley~\cite{CP:1995,CP:1998} in terms of Drinfeld polynomials. We
are interested in the subfamily of such modules which possess a
global crystal basis. Kirillov--Reshetikhin (KR) modules are
finite-dimensional $U'_q(\geh)$-modules $W^{r,s}$ that were
introduced in \cite{HKOTT:2002,HKOTY:1999}. It is expected that each
KR module has a crystal basis $B^{r,s}$, and that every irreducible
finite-dimensional $U'_q(\geh)$-module with crystal basis, is a
tensor product of the crystal bases of KR modules. KR crystals play an
important role in lattice models of statistical mechanics and in the 
Kyoto path construction of highest weight $U_q(\geh)$-modules.

The KR modules $W^{r,s}$ are indexed by a Dynkin node $r$ of the
classical subalgebra (that is, the distinguished simple Lie
subalgebra) $\gc$ of $\geh$ and a positive integer $s$. In general
the existence of $B^{r,s}$ remains an open question. For type
$A_n^{(1)}$ the crystal $B^{r,s}$ is known to
exist~\cite{KKMMNN:1992a} and its combinatorial structure has been
studied~\cite{S:2002}. In many cases, the crystals $B^{1,s}$ and $B^{r,1}$ for
nonexceptional types are also known to exist and their combinatorics has been
worked out in~\cite{KKM:1994,KKMMNN:1992a} and~\cite{JMO:2000,Ko:1999},
respectively. For type $D_n^{(1)}$, $B_n^{(1)}$, and $A_{2n-1}^{(2)}$, which are subject 
of the current paper, the existence of $B^{r,s}$ was recently shown by Okado~\cite{O:2006, OS:2007}.
For the twisted case, this relies on work by Hernandez~\cite{H:2007}.

Viewed as a $U_q(\gc)$-module by restriction, $W^{r,s}$ is generally
reducible. Its decomposition into $U_q(\gc)$-irreducibles was
conjectured in \cite{HKOTT:2002,HKOTY:1999} and verified by
Chari~\cite{C:2001} for the nontwisted cases and Hernandez~\cite{H:2007} for the 
twisted cases.
Kashiwara~\cite{Ka:2002} conjectured that as classical crystals,
many of the KR crystals (the ones conjectured to be perfect in
\cite{HKOTT:2002,HKOTY:1999}) are isomorphic to certain Demazure
subcrystals of affine highest weight crystals. Kashiwara's
conjecture was confirmed by Fourier and Littelmann~\cite{FL:2004} in
the untwisted cases and Naito and Sagaki~\cite{NS:2005} in the
twisted cases. In~\cite{FSS:2006} it was shown that $0$-arrows of these KR crystals
are in fact fixed by the Demazure structure and that this implies that the KR crystals are 
unique if they exist and satisfy certain properties (see Property~\ref{A:KR}).

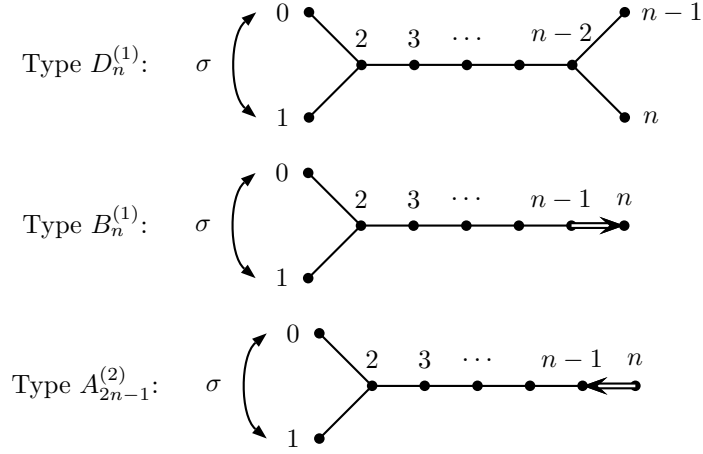
\begin{figure}
\psset{unit=.7cm, dotsize=.2}
\raisebox{0.6cm}{Type $D_n^{(1)}$:} \hspace{1.2cm}
\begin{pspicture}(9,3)
\psdots(1,0)(1,2)(2,1)(3,1)(4,1)(5,1)(6,1)(7,0)(7,2)
\psline(1,0)(2,1)
\psline(1,2)(2,1)
\psline(2,1)(6,1)
\psline(6,1)(7,2)
\psline(6,1)(7,0)
\rput(0.5,2){0}
\rput(0.5,0){1}
\rput(2,1.5){2}
\rput(3,1.5){3}
\rput(4,1.5){$\ldots$}
\rput(5.8,1.6){$n-2$}
\rput(7.9,2){$n-1$}
\rput(7.5,0){$n$}
\pnode(0,0){A}
\pnode(0,2){B}
\ncarc[arcangle=50,nodesep=0mm,arrowsize=.2,arrowinset=0]{<->}{A}{B}
\rput(-1,1){$\sigma$}
\end{pspicture}

\raisebox{0.6cm}{Type $B_n^{(1)}$:} \hspace{1.2cm}
\begin{pspicture}(9,3)
\psdots(1,0)(1,2)(2,1)(3,1)(4,1)(5,1)(6,1)(7,1)
\psline(1,0)(2,1)
\psline(1,2)(2,1)
\psline(2,1)(6,1)
\psline[doubleline=true]{->}(6,1)(7,1)
\rput(0.5,2){0}
\rput(0.5,0){1}
\rput(2,1.5){2}
\rput(3,1.5){3}
\rput(4,1.5){$\ldots$}
\rput(5.8,1.5){$n-1$}
\rput(7,1.5){$n$}
\pnode(0,0){A}
\pnode(0,2){B}
\ncarc[arcangle=50,nodesep=0mm,arrowsize=.2,arrowinset=0]{<->}{A}{B}
\rput(-1,1){$\sigma$}
\end{pspicture}

\raisebox{0.6cm}{Type $A_{2n-1}^{(2)}$:} \hspace{1.2cm}
\begin{pspicture}(9,3)
\psdots(1,0)(1,2)(2,1)(3,1)(4,1)(5,1)(6,1)(7,1)
\psline(1,0)(2,1)
\psline(1,2)(2,1)
\psline(2,1)(6,1)
\psline[doubleline=true]{<-}(6,1)(7,1)
\rput(0.5,2){0}
\rput(0.5,0){1}
\rput(2,1.5){2}
\rput(3,1.5){3}
\rput(4,1.5){$\ldots$}
\rput(5.8,1.5){$n-1$}
\rput(7,1.5){$n$}
\pnode(0,0){A}
\pnode(0,2){B}
\ncarc[arcangle=50,nodesep=0mm,arrowsize=.2,arrowinset=0]{<->}{A}{B}
\rput(-1,1){$\sigma$}
\end{pspicture}
\caption{\label{fig:Dynkin} Dynkin diagram automorphism $\sigma$}
\end{figure}

In this paper, we provide an explicit combinatorial construction of the KR crystals
$B^{r,s}$ of type $D_n^{(1)}$, $B_n^{(1)}$, and $A_{2n-1}^{(2)}$. The construction is
based on the analogue $\sigma$ of the Dynkin diagram automorphism which
interchanges the 0 and the 1 nodes as shown in Figure~\ref{fig:Dynkin}.
A similar construction was given by Shimozono~\cite{S:2002}
for type $A_{n-1}^{(1)}$, for which the Dynkin diagram automorphism maps node
$i$ to node $i+1 \pmod n$. On type $A_{n-1}^{(1)}$ crystals $\sigma$ is the promotion 
operator~\cite{S:2002} and the affine crystal operators are expressed as 
$e_0=\sigma^{-1} \circ e_1 \circ \sigma$ and $f_0=\sigma^{-1} \circ f_1 \circ \sigma$. 
For type $D_n^{(1)}$, the case $B^{2,s}$ was treated in the author's paper with Philip 
Sternberg~\cite{StS:2006}, by providing an explicit description of the Dynkin diagram 
automorphism $\sigma$ that interchanges nodes 0 and 1 on the level of crystals. In his 
thesis~\cite{St:2006}, Sternberg gave a conjecture for the crystal structure for general 
$B^{r,s}$ of type $D_n^{(1)}$. In this paper and~\cite{OS:2007}, we prove Sternberg's conjecture 
by making it more explicit using $\pm$ diagrams and also extend the construction to type 
$B_n^{(1)}$ and $A_{2n-1}^{(2)}$. The $\pm$ diagrams govern the branching of the
underlying classical algebra.

The main result of this paper is the definition of combinatorial crystals $B^{r,s}$
for type $D_n^{(1)}$, $B_n^{(1)}$, and $A_{2n-1}^{(2)}$ given in Definition~\ref{def:comb B}. 
The classical crystal structure is fixed by~\eqref{eq:classical decomp} below and the affine
crystal operators are defined as $e_0=\sigma \circ e_1 \circ \sigma$ and 
$f_0=\sigma \circ f_1 \circ \sigma$ with $\sigma$ as defined in Definition~\ref{def:sigma}. 
As shown in collaboration with Okado~\cite[Theorem 1.2]{OS:2007}, the combinatorial crystals 
$B^{r,s}$ of type $D_n^{(1)}$,  $B_n^{(1)}$, and $A_{2n-1}^{(2)}$ constructed in this paper are in 
fact the Kirillov--Reshetikhin crystals associated with the KR $U'_q(\geh)$-module $W^{r,s}$.
In addition we prove the following theorem. Here $\tilde{B}^{r,s}$ denotes the
unique affine crystal structure of reference~\cite{FSS:2006}  satisfying Property~\ref{A:KR} below. 

\begin{theorem} \label{thm:affine}
For type $D_n^{(1)}$, $B_n^{(1)}$, and $A_{2n-1}^{(2)}$ and $r$ not a spin node,
the crystals $\tilde{B}^{r,s}$ and $B^{r,s}$ are isomorphic.
\end{theorem}

The paper is outlined as follows. In section~\ref{sec:KR} we review the definition of crystals
and the unique characterization of KR crystals coming from Demazure crystal theory as
provided in~\cite{FSS:2006}. In section~\ref{sec:mupad} we briefly describe the implementation
of KR crystals in MuPAD-Combinat; throughout the paper it is demonstrated how to
reproduce examples via MuPAD-Combinat, though the computer implementation is not
used for the proofs. In the section~\ref{sec:construction} the explicit 
construction of $B^{r,s}$ of type $D_n^{(1)}$, $B_n^{(1)}$, and $A_{2n-1}^{(2)}$ is given in 
terms of $\pm$ diagrams. The proof of Theorem~\ref{thm:affine} is provided in section~\ref{sec:proof}.
It is shown in section~\ref{sec:perfect} that $B^{r,s}$ of type $D_n^{(1)}$, $B_n^{(1)}$, and 
$A_{2n-1}^{(2)}$ is perfect. 

\subsection*{Acknowledgements}
I would like to thank Ghislain Fourier, Peter Littelmann, Masato Okado, Mark Shimozono,
and Philip Sternberg for many helpful discussions. Without the collaboration 
on~\cite{StS:2006}, \cite{FSS:2006} and~\cite{OS:2007} the current paper would not have 
been possible. In particular, the notion of $\pm$ diagrams was developed in discussions 
with Mark Shimozono. Many thanks are also due to Christopher Creutzig, 
Fran\c{c}ois Descouens, Teresa Gomez-Diaz, Florent Hivert, and Nicolas Thi\'ery for their support with 
MuPAD-Combinat~\cite{HT:2003}. The implementation of the affine crystal $B^{r,s}$ in
MuPAD-Combinat was essential for the progress of this project!
Thanks to Adrien Boussicault for his help with ps-tricks and Robert Gutierrez for implementing
some algorithms.

\section{Kirillov--Reshetikhin crystals $B^{r,s}$} \label{sec:KR}

Some general definitions regarding crystals are reviewed in sections~\ref{ss:axiom} 
and~\ref{ss:tensor}. The classical crystal structure of type $D_n$, $B_n$, and $C_n$ crystals is 
given in section~\ref{ss:classical}. A unique characterization of Kirillov--Reshetikhin crystals of 
type $D_n^{(1)}$, $B_n^{(1)}$, and $A_{2n-1}^{(2)}$ is reviewed in section~\ref{ss:KR}.

\subsection{Axiomatic definition of crystals} \label{ss:axiom}

Let $\geh$ be a symmetrizable Kac-Moody algebra, $P$ the weight
lattice, $I$ the index set for the vertices of the Dynkin diagram
of $\geh$, $\{\alpha_i\in P \mid i\in I \}$ the simple roots, and
$\{h_i\in P^* \mid i\in I \}$ the simple coroots.
Let $U_q(\geh)$ be the quantized universal enveloping algebra of
$\geh$. A $U_q(\geh)$-crystal \cite{K:1995} is a nonempty set $B$
equipped with maps $\wt:B\rightarrow P$ and
$e_i,f_i:B\rightarrow B\cup\{\emptyset\}$ for all $i\in I$,
satisfying
\begin{align}
\label{eq:e-f}
f_i(b)=b' &\Leftrightarrow e_i(b')=b
\text{ if $b,b'\in B$} \\
\wt(f_i(b))&=\wt(b)-\alpha_i \text{ if $f_i(b)\in B$} \\
\label{eq:string length}
\inner{h_i}{\wt(b)}&=\varphi_i(b)-\varepsilon_i(b).
\end{align}
Here for $b \in B$
\begin{equation*}
\begin{split}
\varepsilon_i(b)&= \max\{n\ge0\mid e_i^n(b)\not=\emptyset \} \\
\varphi_i(b) &= \max\{n\ge0\mid f_i^n(b)\not=\emptyset \}.
\end{split}
\end{equation*}
(It is assumed that $\varphi_i(b),\varepsilon_i(b)<\infty$ for all
$i\in I$ and $b\in B$.) A $U_q(\geh)$-crystal $B$ can be viewed
as a directed edge-colored graph (the crystal graph) whose
vertices are the elements of $B$, with a directed edge from $b$ to
$b'$ labeled $i\in I$, if and only if $f_i(b)=b'$.

Let $\{\Lambda_i \mid i\in I\}$ be the fundamental weights of $\geh$.
For every $b\in B$ define $\varphi(b)=\sum_{i\in I} \varphi_i(b) \La_i$ and
$\varepsilon(b)=\sum_{i\in I} \varepsilon_i(b) \La_i$. An element $b\in B$ is called
highest weight if $e_i(b) = \emptyset$ for all $i\in I$.

The center of an affine Kac--Moody algebra $\mathfrak{g}$ is one-dimensional and 
is generated by the canonical central element $c=\sum_{i\in I}a^\vee_i h_i$, where 
the $a^\vee_i$ are the numbers on the nodes of the Dynkin diagram of the algebra dual 
to $\mathfrak{g}$ given in~\cite[Table Aff of section 4.8]{Kac:1990}.  Moreover, the imaginary
roots of $\mathfrak{g}$ are nonzero integral multiples of the null root
$\delta=\sum_{i\in I} a_i \alpha_i$, where the $a_i$ are the numbers on the nodes of the 
Dynkin diagram of $\mathfrak{g}$ given in~\cite[Table Aff]{Kac:1990}.  Define 
$P_{\cl}=P/\mathbb{Z}\delta$ and 
$P_{\cl}^{+}=\{\la\in P_{\cl} \mid \langle h_{i},\la\rangle \geq 0 \textrm{ for all }i\in I\}$.

\subsection{Tensor products of crystals} \label{ss:tensor}

Let $B_1,B_2,\dotsc,B_L$ be $U_q(\geh)$-crystals. The
Cartesian product $B_L\times \dotsm \times B_2 \times B_1$ has the
structure of a $U_q(\geh)$-crystal using the so-called signature
rule. The resulting crystal is denoted
$B=B_L\otimes\dots\otimes B_2\otimes B_1$ and its elements
$(b_L,\dotsc,b_1)$ are written $b_L\otimes \dotsm \otimes b_1$
where $b_j\in B_j$. The reader is warned that our convention is
opposite to that of Kashiwara \cite{K:1995}. Fix $i\in I$ and
$b=b_L\otimes\dotsm\otimes b_1\in B$. The $i$-signature of $b$ is
the word consisting of the symbols $+$ and $-$ given by
\begin{equation*}
\underset{\text{$\varphi_i(b_L)$ times}}{\underbrace{-\dotsm-}}
\quad \underset{\text{$\varepsilon_i(b_L)$
times}}{\underbrace{+\dotsm+}} \,\dotsm\,
\underset{\text{$\varphi_i(b_1)$ times}}{\underbrace{-\dotsm-}}
\quad \underset{\text{$\varepsilon_i(b_1)$
times}}{\underbrace{+\dotsm+}} .
\end{equation*}
The reduced $i$-signature of $b$ is the subword of the
$i$-signature of $b$, given by the repeated removal of adjacent
symbols $+-$ (in that order); it has the form
\begin{equation*}
\underset{\text{$\varphi$ times}}{\underbrace{-\dotsm-}} \quad
\underset{\text{$\varepsilon$ times}}{\underbrace{+\dotsm+}}.
\end{equation*}
If $\varphi=0$ then $f_i(b)=\emptyset$; otherwise
\begin{equation*}
f_i(b_L\otimes\dotsm\otimes b_1)= b_L\otimes \dotsm \otimes
b_{j+1} \otimes f_i(b_j)\otimes \dots \otimes b_1
\end{equation*}
where the rightmost symbol $-$ in the reduced $i$-signature of
$b$ comes from $b_j$. Similarly, if $\varepsilon=0$ then
$e_i(b)=\emptyset$; otherwise
\begin{equation*}
e_i(b_L\otimes\dotsm\otimes b_1)= b_L\otimes \dotsm \otimes
b_{j+1} \otimes e_i(b_j)\otimes \dots \otimes b_1
\end{equation*}
where the leftmost symbol $+$ in the reduced $i$-signature of $b$
comes from $b_j$. It is not hard to verify that this well-defines
the structure of a $U_q(\geh)$-crystal with
$\varphi_i(b)=\varphi$ and $\varepsilon_i(b)=\varepsilon$ in the above
notation, with weight function
\begin{equation} \label{eq:tensor wt}
\wt(b_L\otimes\dotsm\otimes b_1)=\sum_{j=1}^L \wt(b_j).
\end{equation}
This tensor construction is easily seen to be associative. The
case of two tensor factors is given explicitly by
\begin{equation} \label{eq:f on two factors}
f_i(b_2\otimes b_1) = \begin{cases} f_i(b_2)\otimes b_1
& \text{if $\varepsilon_i(b_2)\ge \varphi_i(b_1)$} \\
b_2\otimes f_i(b_1) & \text{if $\varepsilon_i(b_2)<\varphi_i(b_1)$}
\end{cases}
\end{equation}
and
\begin{equation} \label{eq:e on two factors}
e_i(b_2\otimes b_1) = \begin{cases} e_i(b_2) \otimes b_1 &
\text{if $\varepsilon_i(b_2)>\varphi_i(b_1)$} \\
b_2\otimes e_i(b_1) & \text{if $\varepsilon_i(b_2)\le \varphi_i(b_1)$.}
\end{cases}
\end{equation}

\subsection{Type $D_n$, $B_n$, and $C_n$ crystals} \label{ss:classical}

Crystals of most interest are those associated with a $U_q(\geh)$-module.
In the case when $\geh$ is a simple Lie algebra of nonexceptional type, the crystals 
associated to the $U_q(\geh)$-modules were studied by Kashiwara and Nakashima~\cite{KN:1994}.
Here we review the combinatorial structure in terms of tableaux of the crystals of type $D_n$, $B_n$,
and $C_n$ since these are the finite subalgebras relevant to the KR crystals of type $D_n^{(1)}$, 
$B_n^{(1)}$, and $A_{2n-1}^{(2)}$.

The Dynkin data for type $D_n$, $B_n$, and $C_n$ is given as follows. The simple roots are
\begin{equation*}
\begin{split}
\alpha_i&=\epsilon_i-\epsilon_{i+1}  \qquad \text{for $1\le i<n$}\\
\alpha_n&= \begin{cases} 
		\epsilon_{n-1}+\epsilon_n & \text{for type $D_n$}\\
		\epsilon_n                             &\text{for type $B_n$}\\
		2 \epsilon_n                          & \text{for type $C_n$}
		\end{cases}
\end{split}
\end{equation*}
and the fundamental weights are
\begin{equation*}
\begin{aligned}[3]
\text{Type $D_n$:} \qquad
	&\om_i = \epsilon_1+\cdots+\epsilon_i &&\text{for $1\le i\le n-2$}\\
	&\om_{n-1} = (\epsilon_1+\cdots+\epsilon_{n-1}-\epsilon_n)/2 &&\\
	&\om_n = (\epsilon_1+\cdots+\epsilon_{n-1}+\epsilon_n)/2&&\\[2mm]
\text{Type $B_n$:} \qquad
	&\om_i = \epsilon_1+\cdots+\epsilon_i &&\text{for $1\le i\le n-1$}\\
	&\om_n = (\epsilon_1+\cdots+\epsilon_{n-1}+\epsilon_n)/2&&\\[2mm]
\text{Type $C_n$:} \qquad
	&\om_i = \epsilon_1+\cdots+\epsilon_i &&\text{for $1\le i\le n$}
\end{aligned}
\end{equation*}
where $\epsilon_i\in \mathbb{Z}^n$ is the $i$-th unit standard vector. For type $D_n$ the
nodes $n-1$ and $n$ are spin nodes and for type $B_n$ node $n$ is a spin node. Type $C_n$
does not have any spin node. In fact, $\om_i=\La_i-\langle c,\La_i\rangle \La_0$ are the
level 0 fundamental weights.

Let $X_n=D_n,B_n$, or $C_n$. Any $X_n$ dominant weight $\om$ without a spin component 
can be expressed as $\om=\sum_i c_i \om_i$ for nonnegative integers $c_i$ and the sum runs 
over all $i=1,2,\ldots,n$ not a spin node. In the standard way we represent $\om$ by the partition 
that has exactly $c_i$ columns of height $i$. Conversely, if $\om$ is a partition with at most $n-2$ 
(resp. $n-1$ or $n$) parts we write $c_i(\om)$ for the number of columns of $\om$ of height $i$ for 
$X_n=D_n$ (resp. $X_n=B_n$ or $C_n$). From now on we identify partitions and dominant
weights in this way.

The crystal graph $B(\om_1)$ of the vector representation for type $D_n$, $B_n$, and
$C_n$ is given in Table \ref{tab:vr} by removing the 0 arrows in the crystal $B^{1,1}$
of type $D_n^{(1)}$, $B_n^{(1)}$, and $A_{2n-1}^{(2)}$, respectively.
\begin{table}
\begin{tabular}{|c|l|}
\hline
%
$D_n^{(1)}$ & \raisebox{-1.3cm}{\scalebox{0.7}{
\begin{picture}(365,100)(-10,-50)
\BText(0,0){1} \LongArrow(10,0)(30,0) \BText(40,0){2}
\LongArrow(50,0)(70,0) \Text(85,0)[]{$\cdots$}
\LongArrow(95,0)(115,0) \BText(130,0){n-1}
\LongArrow(143,2)(160,14) \LongArrow(143,-2)(160,-14)
\BText(170,15){n} \BBoxc(170,-15)(13,13)
\Text(170,-15)[]{\footnotesize$\overline{\mbox{n}}$}
\LongArrow(180,14)(197,2) \LongArrow(180,-14)(197,-2)
\BBoxc(215,0)(25,13)
\Text(215,0)[]{\footnotesize$\overline{\mbox{n-1}}$}
\LongArrow(230,0)(250,0) \Text(265,0)[]{$\cdots$}
\LongArrow(275,0)(295,0) \BBoxc(305,0)(13,13)
\Text(305,0)[]{\footnotesize$\overline{\mbox{2}}$}
\LongArrow(315,0)(335,0) \BBoxc(345,0)(13,13)
\Text(345,0)[]{\footnotesize$\overline{\mbox{1}}$}
\LongArrowArc(192.5,-367)(402,69,111)
\LongArrowArcn(152.5,367)(402,-69,-111) \PText(20,2)(0)[b]{1}
\PText(60,2)(0)[b]{2} \PText(105,2)(0)[b]{n-2}
\PText(152,13)(0)[br]{n-1} \PText(152,-9)(0)[tr]{n}
\PText(188,13)(0)[bl]{n} \PText(188,-9)(0)[tl]{n-1}
\PText(240,2)(0)[b]{n-2} \PText(285,2)(0)[b]{2}
\PText(325,2)(0)[b]{1} \PText(192.5,38)(0)[b]{0}
\PText(152.5,-35)(0)[t]{0}
\end{picture}
}}
\\ \hline
%
$B_n^{(1)}$ & \raisebox{-1.3cm}{\scalebox{0.7}{
\begin{picture}(350,100)(-10,-50)
\BText(0,0){1} \LongArrow(10,0)(30,0) \BText(40,0){2}
\LongArrow(50,0)(70,0) \Text(85,0)[]{$\cdots$}
\LongArrow(95,0)(115,0) \BText(125,0){n} \LongArrow(135,0)(155,0)
\BText(165,0){0} \LongArrow(175,0)(195,0) \BBoxc(205,0)(13,13)
\Text(205,0)[]{\footnotesize$\overline{\mbox{n}}$}
\LongArrow(215,0)(235,0) \Text(250,0)[]{$\cdots$}
\LongArrow(260,0)(280,0) \BBoxc(290,0)(13,13)
\Text(290,0)[]{\footnotesize$\overline{\mbox{2}}$}
\LongArrow(300,0)(320,0) \BBoxc(330,0)(13,13)
\Text(330,0)[]{\footnotesize$\overline{\mbox{1}}$}
\LongArrowArc(185,-330)(365,68,112)
\LongArrowArcn(145,330)(365,-68,-112) \PText(20,2)(0)[b]{1}
\PText(60,2)(0)[b]{2} \PText(105,2)(0)[b]{n-1}
\PText(145,2)(0)[b]{n} \PText(185,2)(0)[b]{n}
\PText(225,2)(0)[b]{n-1} \PText(270,2)(0)[b]{2}
\PText(310,2)(0)[b]{1} \PText(185,38)(0)[b]{0}
\PText(145,-35)(0)[t]{0}
\end{picture}
}}
\\ \hline
%
$A_{2n-1}^{(2)}$ & \raisebox{-1.3cm}{\scalebox{0.7}{
\begin{picture}(310,100)(-10,-50)
\BText(0,0){1} \LongArrow(10,0)(30,0) \BText(40,0){2}
\LongArrow(50,0)(70,0) \Text(85,0)[]{$\cdots$}
\LongArrow(95,0)(115,0) \BText(125,0){n} \LongArrow(135,0)(155,0)
\BBoxc(165,0)(13,13)
\Text(165,0)[]{\footnotesize$\overline{\mbox{n}}$}
\LongArrow(175,0)(195,0) \Text(210,0)[]{$\cdots$}
\LongArrow(220,0)(240,0) \BBoxc(250,0)(13,13)
\Text(250,0)[]{\footnotesize$\overline{\mbox{2}}$}
\LongArrow(260,0)(280,0) \BBoxc(290,0)(13,13)
\Text(290,0)[]{\footnotesize$\overline{\mbox{1}}$}
\LongArrowArc(165,-240)(275,65,115)
\LongArrowArcn(125,240)(275,-65,-115) \PText(20,2)(0)[b]{1}
\PText(60,2)(0)[b]{2} \PText(105,2)(0)[b]{n-1}
\PText(145,2)(0)[b]{n} \PText(185,2)(0)[b]{n-1}
\PText(230,2)(0)[b]{2} \PText(270,2)(0)[b]{1}
\PText(165,38)(0)[b]{0} \PText(125,-35)(0)[t]{0}
\end{picture}
}}
\\ \hline
\end{tabular}\vspace{4mm}
\caption{\label{tab:vr}KR crystal $B^{1,1}$}
\end{table}
The crystal $B(\om_\ell)$ for $\ell$ not a spin node can be realized as the connected component
of $B(\om_1)^{\otimes \ell}$ containing the element $\ell\otimes \ell-1\otimes \cdots \otimes 1$.
Similarly, the crystal $B(\om)$ labeled by a dominant weight $\om=\om_{\ell_1}+\cdots+\om_{\ell_k}$
with $\ell_1\ge \ell_2 \ge \cdots \ge \ell_k$ not containing spin nodes can be realized as the connected 
component in $B(\om_{\ell_1}) \otimes \cdots \otimes B(\om_{\ell_k})$ containing the element 
$u_{\om_{\ell_1}} \otimes \cdots \otimes u_{\om_{\ell_k}}$, where $u_{\om_i}$ is the highest 
weight element in $B(\om_i)$. As shown in~\cite{KN:1994},  the elements of $B(\om)$ can be labeled 
by tableaux of shape $\om$ in the alphabet $\{1,2,\ldots, n, \overline{n}, \ldots, \overline{1}\}$
for types $D_n$ and $C_n$ and the alphabet $\{1,2,\ldots, n, 0, \overline{n}, \ldots, \overline{1}\}$
for type $B_n$. For the explicit rules of type $D_n$, $B_n$, and $C_n$ tableaux we refer the reader 
to~\cite{KN:1994}; see also~\cite{HongKang:2002}.

\subsection{Kirillov--Reshetikhin crystal of type $D_n^{(1)}$, $B_n^{(1)}$, and $A_{2n-1}^{(2)}$} 
\label{ss:KR}

Let $\geh$ be an affine Kac--Moody algebra of rank $n$.
Kirillov--Reshetikhin (KR) modules $W^{r,s}$, labeled by a positive integer $s$ and 
$r \in \{1,2,\ldots,n\}$, are finite-dimensional $U'_q(\geh)$-modules. See~\cite{HKOTT:2002} for the
precise definition. It was shown in collaboration with Okado~\cite{O:2006, OS:2007} that $W^{r,s}$ 
of type $D_n^{(1)}$, $B_n^{(1)}$, and $A_{2n-1}^{(2)}$ has a global crystal basis $B^{r,s}$. 

In~\cite{HKOTT:2002} a conjecture is given for the decomposition of each KR module $W^{r,s}$ 
into its $\geh_0$-components, which was proven by Chari~\cite{C:2001} for the nonexceptional
untwisted algebras and Hernandez~\cite{H:2007} for the twisted cases.
For $\geh$ of type $D_n^{(1)}$, $B_n^{(1)}$, or $A_{2n-1}^{(2)}$
the underlying finite Lie algebra $\geh_0$ is of type $X_n=D_n,B_n$, or $C_n$, respectively.
Explicitly, as a $X_n$-crystal, the KR crystal $B^{r,s}$ of type $D_n^{(1)}$, $B_n^{(1)}$ or 
$A_{2n-1}^{(2)}$ decomposes into the following irreducible components
\begin{equation} \label{eq:classical decomp}
	B^{r,s} \cong \bigoplus_{\om} B(\om),
\end{equation}
for $1\le r\le n$ not a spin node. Here $B(\om)$ is the $X_n$-crystal of highest weight $\om$
and the sum runs over all dominant weights $\om$ that can be obtained from
$s\om_r$ by the removal of vertical dominoes, where $\om_i$ are the fundamental weights of 
$X_n$ as defined in section~\ref{ss:classical}. 

In~\cite[Corollary 4.6]{FSS:2006} we showed that Property~\ref{A:KR} below defines a unique
affine crystal structure $\tilde{B}^{r,s}$. In this paper we construct the explicit combinatorial
affine crystal structure $B^{r,s}$ of $\tilde{B}^{r,s}$ for type $D_n^{(1)}$,
$B_n^{(1)}$, and $A_{2n-1}^{(2)}$ by showing that it satisfies these properties.
In a subsequent paper with Okado~\cite[Theorem 1.2]{OS:2007} it is shown that the 
combinatorial crystals $B^{r,s}$ of this paper are indeed the crystals associated to the KR 
modules $W^{r,s}$. Combining these two results implies in particular that the KR crystals
$B^{r,s}$ satisfy the conditions of Property~\ref{A:KR}.

Before stating the Property~\ref{A:KR}, we need to make a few definitions.
An $I$-crystal $B$ is \textit{regular} if, for each subset $K\subset
I$ with $|K|=2$, each $K$-component of $B$ is isomorphic to the
crystal basis of an irreducible integrable highest weight
$U'_q(\geh_K)$-module where $\geh_K$ is the subalgebra of $\geh$
with simple roots $\alpha_i$ for $i\in K$.

The Dynkin diagram of type $D_n^{(1)}$, $B_n^{(1)}$, and $A_{2n-1}^{(2)}$ all have an 
automorphism $\sigma$ interchanging nodes 0 and 1. See Figure~\ref{fig:Dynkin}.

\begin{property} \label{A:KR} Let $\tilde{B}^{r,s}$ of type $D_n^{(1)}$, $B_n^{(1)}$, and
$A_{2n-1}^{(2)}$ be the crystal with the following properties:
\begin{enumerate}
\item \label{A:classical} As an $X_n$ crystal, $\tilde{B}^{r,s}$ decomposes according 
to~\eqref{eq:classical decomp}, where $X_n=D_n,B_n$, and $C_n$, respectively.
\item \label{A:regular} $\tilde{B}^{r,s}$ is regular.
\item \label{A:u} There is a unique element $u\in \tilde{B}^{r,s}$ such that
\begin{align}\label{eq:u}
\ve(u)=s\La_0 \quad \text{and} \quad 
\vp(u) = \begin{cases} s \La_0 & \text{for $r$ even},\\ s \La_1 & \text{for $r$ odd.} \end{cases}
\end{align}
\item \label{A:auto} $\tilde{B}^{r,s}$ admits the automorphism corresponding to $\sigma$ 
(also denoted $\sigma$) such that 
$\ve\circ \sigma = \sigma\circ \ve$ and $\vp\circ\sigma=\sigma\circ\vp$.
\end{enumerate}
\end{property}

\begin{theorem} \cite[Corollary 4.6]{FSS:2006} \label{thm:specification}
The crystal $\tilde{B}^{r,s}$ of type $D_n^{(1)}$, $B_n^{(1)}$, and $A_{2n-1}^{(2)}$ is uniquely 
determined by the conditions of Property~\ref{A:KR}.
\end{theorem}

\section{MuPAD-Combinat implementation} \label{sec:mupad}

\begin{figure}
\scalebox{.5}{\includegraphics{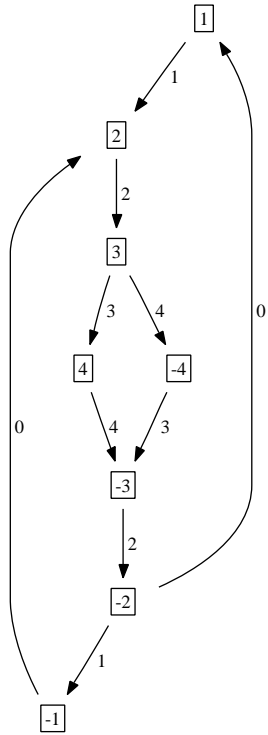}}
\caption{Crystal graph $B^{1,1}$ of type $D_4^{(1)}$ obtained via MuPAD-Combinat
\label{fig:crystal graph}}
\end{figure}

The Kirillov--Reshetikhin crystals $B^{r,s}$ of type $A_n^{(1)}$, $B_n^{(1)}$, $D_n^{(1)}$, and 
$A_{2n-1}^{(2)}$ are implemented in MuPAD-Combinat~\cite{HT:2003}, 
an open source algebraic combinatorics package for the computer algebra
system MuPAD~\cite{Fuchssteiner:1996}. A KR crystal is declared via the command
\begin{center}
	combinat::crystals::kirillovReshetikhin($r$,$s$,type)
\end{center}
For example
\begin{verbatim}
	>> KR:=combinat::crystals::kirillovReshetikhin(2,2,["D",4,1]):
\end{verbatim}
defines the Kirillov--Reshetikhin crystal $B^{2,2}$ of type $D_4^{(1)}$.
An element $t$ in this crystal is specified by its corresponding tableau via
\begin{verbatim}
>> t:=KR([[3],[1]])

                                     +---+
                                     | 3 |
                                     +---+
                                     | 1 |
                                     +---+
\end{verbatim}
Then $e_0(t)$ and $\sigma(t)$ are obtained via
\begin{verbatim}
>> t::e(0)

                                     +----+
                                     | -2 |
                                     +----+
                                     | 3  |
                                     +----+
\end{verbatim}
and
\begin{verbatim}
>> t::sigma()

                                  +----+----+
                                  | -2 | -1 |
                                  +----+----+
                                  | 2  | 3  |
                                  +----+----+
\end{verbatim}
The full crystal graph can be obtained by
\begin{verbatim}
>> KR:=combinat::crystals::kirillovReshetikhin(1,1,["D",4,1]):
>> KR::crystalGraph("filename.dot"):
\end{verbatim}
followed by {\tt dot -Tpdf filename.dot -o filename.pdf} in the command line.
The result is presented in Figure~\ref{fig:crystal graph}~\cite{graphviz} .
Throughout the paper more functionalities of the implementation are presented by
examples. It should be noted that the computer implementation is not used for the proofs,
but rather for illustration and intuition.

\section{Explicit construction of $B^{r,s}$} \label{sec:construction}

In this section we define a combinatorial crystal $B^{r,s}$ which satisfies all conditions
of Property~\ref{A:KR}.
By~\eqref{eq:classical decomp}, the $X_n$-crystal structure of $B^{r,s}$ is fixed, that is, all Kashiwara
operators $e_i$ and $f_i$ for $1\le i\le n$ are determined. Hence the complete crystal
structure of $B^{r,s}$ is determined by specifying the affine Kashiwara crystal operators
$e_0$ and $f_0$. We define $e_0$ and $f_0$ by constructing the crystal analogue
$\sigma$ of the automorphism of the Dynkin diagram that interchanges the 0 and 1 node (see 
Definition~\ref{def:sigma} below). Then
\begin{equation} \label{eq:e0}
\begin{split}
f_0 &= \sigma \circ f_1 \circ \sigma,\\
e_0 &= \sigma \circ e_1 \circ \sigma.
\end{split}
\end{equation}
\begin{definition} \label{def:comb B}
The combinatorial crystal $B^{r,s}$ is given by the classical 
decomposition~\eqref{eq:classical decomp} and the affine crystal operators
$e_0=\sigma\circ e_1\circ \sigma$ and $f_0=\sigma\circ f_1\circ \sigma$ with
$\sigma$ as defined in Definition~\ref{def:sigma}.
\end{definition}

The involution $\sigma$ is first defined on $X_{n-1}$ highest weight elements and then extended
to any element in $B^{r,s}$ (see section~\ref{ss:sigma}). In section~\ref{ss:branching} we discuss 
how the branching from $X_n$ to $X_{n-1}$ can be formulated in terms of $\pm$ diagrams.  
In section~\ref{ss:bij} a bijection between $X_{n-1}$ highest weight elements and $\pm$ diagrams 
is given.

\subsection{$X_n$ to $X_{n-1}$ branching}
\label{ss:branching}

We introduce combinatorial objects called $\pm$ diagrams to describe
the branching from $X_n$ to the subalgebra of type $X_{n-1}$
obtained by removing the Dynkin node $1$. A $\pm$ diagram $P$ of
shape $\La/\la$ is a sequence of partitions $\la\subset \mu \subset \La$ 
such that $\La/\mu$ and $\mu/\la$ are horizontal strips. We
depict this $\pm$ diagram by the skew tableau of shape $\La/\la$ in
which the cells of $\mu/\la$ are filled with the symbol $+$ and
those of $\La/\mu$ are filled with the symbol $-$. Write
$\La=\os(P)$ and $\la=\is(P)$ for the outer and inner shapes of the
$\pm$ diagram $P$.

If $\la$ is a dominant weight for the simple Lie algebra $\geh$
write $B_\geh(\la)$ for the crystal graph of the highest weight
$U_q(\geh)$-module of highest weight $\la$.

\begin{proposition} \label{P:branch}
For an $X_n$ dominant weight $\La$ with no spin weights, there is an
isomorphism of $X_{n-1}$-crystals
\begin{align*}
  B_{X_n}(\La) \cong \bigoplus_{\substack{\text{$\pm$ diagrams $P$} \\ \os(P)=\La}}
  B_{X_{n-1}}(\is(P)).
\end{align*}
That is, the multiplicity of $B_{X_{n-1}}(\la)$ in $B_{X_n}(\La)$,
is the number of $\pm$ diagrams of shape $\La/\la$.
\end{proposition}

\begin{proof}
This follows directly from the branching rules for $X_n$ to $X_{n-1}$
(see for example~\cite[pg. 426]{FH:1991}).
\end{proof}

\begin{example}
Let $\La=\om_4+\om_2$. The corresponding $\pm$ diagrams are
\begin{equation*}
\begin{split}
&\yng(1,1,2,2) \quad \young(+,\mb,\mb\mb,\mb\mb) \quad
\young(\mb,\mb,\mb+,\mb\mb) \quad \young(+,\mb,\mb+,\mb\mb) \quad
\young(-,\mb,\mb\mb,\mb\mb) \quad
\young(\mb,\mb,\mb-,\mb\mb) \quad \young(-,\mb,\mb-,\mb\mb) \quad
\young(-,+,\mb\mb,\mb\mb) \quad \young(+,\mb,\mb-,\mb\mb) \quad\\
& \young(-,+,\mb-,\mb\mb) \quad \young(-,\mb,\mb+,\mb\mb) \quad
\young(\mb,\mb,\mb-,\mb+) \quad
\young(-,\mb,\mb-,\mb+) \quad \young(-,+,\mb+,\mb\mb) \quad
\young(+,\mb,\mb-,\mb+) \quad \young(-,+,\mb-,\mb+)
\end{split}
\end{equation*}
Therefore (suppressing the subscript $X_{n-1}$) we have
\begin{align*}
  B_{X_n}(\om_4+\om_2) &\cong
  B(\om_4+\om_2) \oplus B(\om_4+\om_1)^{\oplus 2} \oplus  B(\om_4)
  \oplus B(\om_3+\om_2)^{\oplus 2} \oplus \\
  &\quad \,\,B(\om_3+\om_1)^{\oplus 4} \oplus B(\om_3)^{\oplus 2}
  \oplus B(2\om_2) \oplus B(\om_2+\om_1)^{\oplus 2} \oplus B(\om_2).
\end{align*}
\end{example}

\subsection{$X_{n-1}$ highest weight elements} \label{ss:bij}
We now suppose that $B^{r,s}$ is a KR crystal of type $D_n^{(1)}$, $B_n^{(1)}$, or
$A_{2n-1}^{(2)}$ with $1\le r \le n$ not a spin node and $\La$ is such that $B_{X_n}(\La)$ is a
$X_n$-summand in~\eqref{eq:classical decomp}, where $X_n=D_n,B_n$, or $C_n$, respectively.
That is, $\La$ has exactly $s$ columns, and for each column, its height is at most $r$ and of the 
same parity as $r$. If $r$ is even then columns of height zero are
allowed.

Let $\la$ be a partition associated with a ``spinless"
$X_{n-1}$-highest weight. We now give an explicit bijection
$\bij:P\mapsto b$ from the set of $\pm$ diagrams $P$ of shape $\La/\la$
to the set of $X_{n-1}$-highest weight vectors $b$ of
$X_{n-1}$-weight $\la$ in $B_{X_n}(\La)$ (that is, the elements
$b\in B_{X_n}(\La)$ such that $\ve_i(b)=0$ and $\vp_i(b)=c_i(\la)$ 
for all $2\le i\le n$).
This bijection has the property of preserving the $1$-weight, that
is, $m_+(P)-m_-(P)=m_1(b)-m_{\ba{1}}(b)$ where $m_+(P)$ and $m_-(P)$
are respectively the number of symbols $+$ and $-$ in $P$. Here, for
$b\in B_{X_n}(\La)$ let $m_i(b)$ (resp. $m_{\ba{i}}(b)$) be the
number of symbols $i$ (resp. $\ba{i}$) that occur in $b$, for $1\le
i\le n$. The bijection also satisfies the property that the barred letters 
in $b$ occur precisely in the positions in $P$ containing the symbol $-$. 

The precise correspondence is given by the following algorithm.
Start with a $\pm$ diagram $P$ of shape $\La/\la$. Let
$d=m_+(T)-m_-(T)$. We build the tableau $b\in B_{X_n}(\La)$ as
follows. Place the maximum number of symbols $1$ and $\bar{1}$ into
the diagram of $\La$, placing the $1$s into the first row and the
$\bar{1}$s into the rightmost positions of $-$ in $P$, such that
$m_1(b)-m_{\ba{1}}(b)=d$, and the placement is legal in a
$X_n$-tableau. Assume that all letters
$1,\ba{1},2,\ba{2},\dotsc,i-1,\ba{i-1}$ have already been placed.
Fill the remainder of row $i-1$ with $i$'s. Then place the maximum
number of symbols $i$ into row $i$ and symbols $\bar{i}$ into the
available $-$ positions starting from the right, such that the
$i$-weight is correct ($m_i(b)-m_{\ba{i}}(b)=\la_{i-1}$) and the
partial tableau is legal in a $X_n$-tableau.

\begin{remark} During the placement of the letters $i$ and $\ba{i}$
in the computation of the map $P\mapsto b$ it suffices to check the
following conditions:
\begin{enumerate}
\item The unbarred subtableau is semistandard with respect to the
total order $1<2<\dotsm<i$;
\item A column containing $\bar{i}$ cannot contain all of the
letters
$1,2,\dotsc,i$;
\item there are no configurations of the form
\begin{equation*}
\begin{array}{c|c} & \bar{i}\\ &\\ i&i\end{array} \qquad
\text{or} \qquad
\begin{array}{c|c} \bar{i}&\bar{i}\\ & \\ i& \end{array}
\end{equation*}
in adjacent columns.
\end{enumerate}
\end{remark}

Alternatively, the bijection $\bij:P\mapsto b$ from $\pm$ diagrams $P$ of shape 
$\La/\la$ to the set of $X_{n-1}$-highest weight vectors $b$ of
$X_{n-1}$-weight $\la$ is as follows. Place $\ba{1}$ in all positions in $P$
that contain a $-$ and fill the remainder of all columns by strings of the 
form $23\ldots k$. We move through the columns of $b$  from top to bottom,
left to right. Each $+$ in $P$ (starting with the leftmost moving to the right)
will alter $b$ as we move through the columns. Suppose the $+$ is at height $h$ in $P$.
If one encounters a $\ba{1}$, replace $\ba{1}$ by $\ba{h+1}$. If one encounters a $2$,
replace the string $23\ldots k$ by $12\ldots h h+2\ldots k$.

\begin{example}
The $X_{n-1}$-highest weight elements corresponding to
\begin{equation*}
\young(+,\mb,\mb+,\mb\mb) \qquad
\young(+,\mb,\mb-,\mb\mb) \qquad
\young(-,\mb,\mb+,\mb\mb) \qquad
\young(-,\mb,\mb-,\mb\mb)
\end{equation*}
are
\begin{equation*}
\young(4,3,22,11) \qquad
\young(4,3,2\ab,12) \qquad
\young(\cb,4,33,22) \qquad
\young(\ab,4,3\ab,22)
\end{equation*}
respectively.
\end{example}

\begin{example} \label{ex:hw}
The $X_{n-1}$-highest weight element corresponding to
\begin{equation*}
\young(+-,\mb+,\mb\mb--,\mb\mb\mb+) \quad \text{under $\bij$ is} \quad
\young(4\db,34,23\ab\ab,1122)\,.
\end{equation*}
In MuPAD-Combinat this example can be reproduced by
\begin{verbatim}
>> KR:=crystals::kirillovReshetikhin(4,5,["D",6,1]):
>> P:=[["+","-"],["","+"],["","","-","-"],["","","","+"]]:
>> KR::diagramsToTableaux(P)
\end{verbatim}
inside the KR crystals $B^{4,5}$ of type $D_6^{(1)}$.
\end{example}

For a given $\pm$ diagram $P$ one can easily construct the string of 
operators $f_{\aar} := f_{a_1} f_{a_2} \cdots f_{a_\ell}$ such that
$\bij(P) = f_{\aar} u$, where $u$ is the highest weight vector in $B_{X_n}(\Lambda)$. 
Start with the empty string $\aar=()$. Scan the columns of $P$ from right
to left. For each column of $P$ for which a $+$ can be added, append
$(1,2, \ldots, h)$ to $\aar$, where $h$ is the height of the added $+$.
Next scan $P$ from left to right and for each column that contains a $-$ in $P$,
append to $\aar$ the string $(1,2,\ldots,n,n-2,n-3,\ldots, h)$ for type $D_n$,
$(1,2,\ldots,n-1,n,n,n-1,\ldots,h)$ for type $B_n$, and $(1,2,\ldots,n-1,n,n-1,\ldots,h)$ for
type $C_n$, where $h$ is the height of the $-$ in $P$.
 
\begin{proposition} \label{prop:to hw}
Let $P$ be a $\pm$ diagram of shape $\Lambda/\lambda$. Let $\aar$ be as constructed
above. Then $\bij(P) = f_{\aar} u$, where $u$ is the $X_n$ highest weight vector in $B_{X_n}(\La)$. 
\end{proposition}
\begin{proof}
This follows immediately from the previous discussion.
\end{proof}

\begin{example} \label{ex:string}
Let $P$ be the $\pm$ diagram of Example~\ref{ex:hw}. Then for type $D_6$ we have
$\aar=(1,\hspace{0.2cm}1,2,3,4,5,6,4,\hspace{0.2cm}1,2,3,4,5,6,4,3,2,\hspace{.2cm} 
1,2,3,4,5,6,4,3,2)$.
The command to obtain the string $\aar$ in MuPAD-Combinat is
\begin{verbatim}
>> KR::diagramsString(P,6) 
\end{verbatim}
\end{example}

\subsection{Definition of $\sigma$} \label{ss:sigma}
Let $\sigma$ be the automorphism of the Dynkin diagram that interchanges 0 and 1 as shown
in Figure~\ref{fig:Dynkin}. Here we construct the analogous automorphism
on $B^{r,s}$ which, by abuse of notation, is also called $\sigma$. By construction $\sigma$ is
an involution which commutes with $f_i$ and $e_i$ for $i=2,3,\ldots,n$. Hence it
suffices to define $\sigma$ on $X_{n-1}$ highest weight elements
where $X_{n-1}$ is the subalgebra whose Dynkin diagram is obtained
from that of $X_n$ by removing node $1$. Because of the
bijection $\bij$ between $\pm$ diagrams and $X_{n-1}$-highest weight
elements defined in section~\ref{ss:bij}, it suffices to define the map on $\pm$ diagrams.

Let $P$ be a $\pm$ diagram of shape $\La/\la$. Let $c_i=c_i(\la)$ be
the number of columns of height $i$ in $\la$ for all $1\le i<r$ with
$c_0=s-\la_1$. If $i\equiv r-1 \pmod{2}$, then in $P$, above each
column of $\la$ of height $i$, there must be a $+$ or a $-$. Suppose there are
$p_i$ of the $+$ symbols above the columns of height $i$, so that there are $(c_i-p_i)$ of
the $-$ symbols. Change this to $(c_i-p_i)$ of the $+$ symbols and $p_i$ of the $-$ symbols.
If $i\equiv r \pmod{2}$, then in $P$, above each column of $\la$ of height $i$,
either there are no signs or a $\mp$ pair. Suppose there are $p_i$
$\mp$ pairs above the columns of height $i$. Change this to
$(c_i-p_i)$ $\mp$ pairs. The result is $\sigD(P)$, which has the
same inner shape $\la$ as $P$ but a possibly different outer shape.

\begin{example}\label{ex:pm}
Let $s=5$. Then we have
\begin{equation*}
P=\young(+-,\mb+,\mb\mb--,\mb\mb\mb+)\qquad
\sigD(P)=\young(-,\mb,\mb\mb+-,\mb\mb\mb+)\,.
\end{equation*}
In MuPAD-Combinat this is achieved via
\begin{verbatim}
>> KR::sigmaOnDiagrams(P)
\end{verbatim}
\end{example}

To define $\sigma$ on any element $b\in B^{r,s}$, first go to the $X_{n-1}$-highest
weight element in the component of $b$ using the crystal raising operators $e_i$, then 
use $\sigD$ on the corresponding $\pm$-diagram, and go back to the same element within the new
$X_{n-1}$ component using the string of $f_i$ that reverses the previously applied string of $e_i$.

\begin{definition} \label{def:sigma}
Let $b\in B^{r,s}$ and $e_{\aar} := e_{a_1} e_{a_2} \cdots e_{a_\ell}$ be such that
$e_{\aar}(b)$ is a $X_{n-1}$ highest weight crystal element. Define 
$f_{\aal}:= f_{a_\ell} f_{a_{\ell-1}} \cdots f_{a_1}$. Then 
\begin{equation} \label{eq:def sigma}
\sigma(b) := f_{\aal} \circ \bij \circ \sigD \circ \bij^{-1} \circ e_{\aar}(b),
\end{equation}
with $\sigD$ as defined above and $\bij$ as defined in section~\ref{ss:bij}.
\end{definition}

\begin{example} Take
\begin{equation*}
b=\young(\db\bb,34,23\ab\ab,1123) \; \in \; B^{4,5}
\end{equation*}
in type $D_6^{(1)}$. Then the corresponding $D_5$-highest weight vector
$e_{\aar}(b)$ is the tableau of Example~\ref{ex:hw} where $\aar=(4,6,5,4,3,2,2)$.
Using the results of Examples~\ref{ex:hw} and~\ref{ex:pm} we obtain
\begin{equation*}
\sigma(b) = \young(\bb,\db,334\ab,1223)\; .
\end{equation*}
In MuPAD-Combinat this example can be checked by
\begin{verbatim}
>> b:=KR([[-4,-2],[3,4],[2,3,-1,-1],[1,1,2,3]]):
>> b::sigma()
\end{verbatim}
\end{example}

Note that the shape of $b\in B^{r,s}$, or equivalently, the classical component $B(\om)$
in~\eqref{eq:classical decomp} will in general change under the application of $\sigma$.

\section{Proof of Theorem~\ref{thm:affine}} \label{sec:proof}

By Theorem~\ref{thm:specification} we know that $B^{r,s}$ is uniquely determined by the
four conditions of  Property~\ref{A:KR}. Hence to prove Theorem~\ref{thm:affine} it suffices 
to show that the crystal given by~\eqref{eq:classical decomp} together with the affine crystal 
operators $e_0$ and $f_0$ of~\eqref{eq:e0} with $\sigma$ as in Definition~\ref{def:sigma} 
satisfies the conditions~(\ref{A:classical}), (\ref{A:regular}), (\ref{A:u}) and~(\ref{A:auto}) of 
Property~\ref{A:KR}.

By construction Property~\ref{A:KR} (\ref{A:classical}) holds and $\sigma$ satisfies 
Property~\ref{A:KR} (\ref{A:auto}).

If $r$ is even, let $u$ be the unique element in $B(\emptyset)$. If $r$ is odd, let
$u$ be the $X_n$-highest weight element in the component $B(s\om_1)$.
It is not hard to check that $u$ is the element of Property~\ref{A:KR} (\ref{A:u}).
Because of the condition $\ve(u)=s\Lambda_0$, $u$ has to be a classical
highest weight vector. Each classical highest weight vector $v\in B(\La)$
corresponds to the $\pm$ diagram of outer shape $\La$ with a $+$ in every column.
It can be shown explicitly, that $u$ is the unique highest weight vector $v$
satisfying $\vp(u)=s\La_0$ for $r$ even and $s\La_1$ for $r$ odd.

It remains to consider the regularity Property~\ref{A:KR} (\ref{A:regular}).
(In principle this follows from~\cite[Theorem 1.2]{OS:2007} which shows that 
$B^{r,s}$ is the crystal corresponding to the KR module $W^{r,s}$; the proof
in~\cite{OS:2007} uses only a special case of Lemma~\ref{lem:e1 action} below. We give
an independent proof here).
Since $B^{r,s}$ decomposes as~\eqref{eq:classical decomp} as an $X_n$-crystal, each 
$K$-component of $B^{r,s}$ is an irreducible integrable highest weight $U_q(\geh_K)$-crystal 
for $K=\{i,j\}$ with $i,j\neq 0$. For $K=\{0,i\}$ with $i\neq 1$ this is also true since
\begin{equation*}
e_0 e_i = \sigma e_1\sigma e_i = \sigma(e_1 \sigma e_i \sigma ) \sigma
= \sigma (e_1 e_i) \sigma,
\end{equation*}
so that by conjugation with $\sigma$ the crystal behaves just like a $\{1,i\}$-crystal.
Hence it suffices to consider the case $K=\{0,1\}$ in which case we need to show that
$e_0$ and $e_1$ commute since the 0 and 1 node are not adjacent in the
$D_n^{(1)}$, $B_n^{(1)}$, and $A_{2n-1}^{(2)}$ Dynkin diagrams.

For us $I=\{0,1,2,3,\ldots,n\}$ are the labels of the affine Dynkin diagram. In the following we 
will also consider the subalgebras $X_n$ with Dynkin set $I_{X_n} = \{1,2,3,\ldots,n\}$, 
$X_{n-1}$ with Dynkin set $I_{X_{n-1}} = \{2,3,\ldots,n\}$, and $X_{n-2}$ with Dynkin set 
$I_{X_{n-2}} = \{3,4,\ldots,n\}$. Since $e_i$ with $i\in I_{X_{n-2}}$ commutes with both $e_0$ 
and $e_1$, it suffices to prove the commutativity of $e_0$ and $e_1$ for $X_{n-2}$ highest 
weight vectors. As we have seen in section~\ref{ss:bij}, the $X_{n-1}$-highest weight elements 
in the branching $X_n\to X_{n-1}$ can be described by $\pm$ diagrams. Similarly the 
$X_{n-2}$-highest weight elements in the branching $X_{n-1}\to X_{n-2}$ can be described 
by $\pm$ diagrams. Hence each $X_{n-2}$-highest weight vector is uniquely determined by a 
pair of $\pm$ diagrams $(P,p)$ such that $\is(P)=\os(p)$. The diagram $P$ specifies the 
$X_{n-1}$-component $B_{X_{n-1}}(\is(P))$ in $B_{X_n}(\os(P))$, and $p$ specifies the 
$X_{n-2}$ component inside $B_{X_{n-1}}(\is(P))$. Let $\Psi$ denote the map $(P,p) \mapsto b$ 
from a pair of $\pm$ diagrams to a $X_{n-2}$ highest weight vector.

\begin{example} \label{ex:Psi}
The tableau $b\in B^{4,3}$ of type $D_6^{(1)}$ corresponding to 
\begin{equation*}
P=\young(-,+,\mb+-,\mb\mb\mb) \quad \text{and} \quad
p=\young(-,\mb\mb+) \quad \text{is} \quad
\Psi(P,p)=\young(\cb,\db,34\ab,133) \;.
\end{equation*}
In MuPAD-Combinat this is achieved via
\begin{verbatim}
>> KR:=crystals::kirillovReshetikhin(4,3,["D",6,1]):
>> P:=[["-"],["+"],["","+","-"],["","",""]]:
>> p:=[["-"],["","","+"]]:
>> KR::twoDiagramsToTableaux(P,p)
\end{verbatim}
\end{example}

By \eqref{eq:e0} the commutation of $e_0$ and $e_1$ is equivalent to the commutativity
of the following diagram
\begin{equation}\label{eq:e1 sigma}
\xymatrix{
 B^{r,s} \ar[d]_{\sigma} \ar[r]^{e_1} & B^{r,s} \cup \{\emptyset\} \ar[r]^{\sigma} & 
 B^{r,s} \cup \{\emptyset\} \ar[r]^{e_1} & B^{r,s} \cup \{\emptyset\} \ar[d]^{\sigma}\\
 B^{r,s} \ar[r]_{e_1} & B^{r,s} \cup \{\emptyset\} \ar[r]_{\sigma} & B^{r,s}  \cup \{\emptyset\}
  \ar[r]_{e_1} & B^{r,s} \cup \{\emptyset\}
}
\end{equation}

By~\eqref{eq:def sigma}, the action of $\sigma$ on a pair of $\pm$ diagrams $(P,p)$ is given by
$(\sigD(P),p)$. 

The operator $e_1$ will either change a 2 into a 1 or a $\ab$ into a $\bb$
in $b=\Psi(P,p)$. On the level of $(P,p)$ this means that either a $+$ from $p$ transfers to $P$,
or a $-$ moves from $P$ to $p$. To describe the precise action of $e_1$ on $(P,p)$ perform
the following algorithm:
\begin{enumerate}
\item Successively run through all $+$ in $p$ from left to right and, if possible, pair it with 
the leftmost yet unpaired $+$ in $P$ weakly to the left of it.
\item Successively run through all $-$ in $p$ from left to right and, if possible, pair it with
the rightmost  yet unpaired $-$ in $P$ weakly to the left.
\item Successively run through all yet unpaired $+$ in $p$ from left to right and, if possible,
pair it with the leftmost yet unpaired $-$ in $p$.
\end{enumerate}

\begin{lemma} \label{lem:e1 action}
If there is an unpaired $+$ in $p$,  $e_1$ moves the rightmost unpaired $+$ in $p$ to $P$. 
Else, if there is an unpaired $-$ in $P$, $e_1$ moves the leftmost unpaired $-$ in $P$ to $p$.
Else $e_1$ annihilates $(P,p)$.
\end{lemma}

\begin{proof}[Sketch of proof]
Let $\aar$ be the string of indices in $I_{X_n}=\{1,2,\ldots, n\}$ of Proposition~\ref{prop:to hw} such 
that $\bij(P)=f_{\aar} u$, where $u$ is the $X_n$-highest weight vector of $B_{X_n}(\os(P))$.
Similarly, let $\bbr$ be the string of indices in $I_{X_{n-1}}=\{2,3,\ldots,n\}$ such that
$\Psi(P,p) = f_{\bbr} f_{\aar} u$. Again $\bbr$ is determined by Proposition~\ref{prop:to hw},
where this time $f_{\aar}u$ takes the role of the $X_{n-1}$-highest weight vector of 
$B_{X_{n-1}}(\is(P))$.

In the following we give the arguments for $X_n=D_n$. The cases $X_n=B_n,C_n$ can be
treated in the same fashion.
Recall that the action of $f_i$ on a tableau can be determined by reading the tableau columnwise.
After determining the reduced $i$-signature, $f_i$ acts on the rightmost possible letter;
see Sections~\ref{ss:tensor} and~\ref{ss:classical}.
Match each string $2,3,\ldots,n,n-2,\ldots,h'$ in $\bbr$ corresponding to a $-$ in $p$ 
with a string $1,2,\ldots,n,n-2,\ldots,h$ in $\aar$ which corresponds to a $-$ in a column in $P$
weakly to the left, meaning that $h'<h$. The raising operator $e_1$ cannot act on such a 
combination. This corresponds to the pairing of step (2). Similarly, match each string $2,3,\ldots,h'$ 
in $\bbr$ with a string $1,2,\ldots,h$ in $\aar$, where $h'> h$. The operator $e_1$ cannot act on such
a combination. This is equivalent to the pairing of step (1) (recall that the strings $(1,)2,3,\ldots,h$
correspond to columns in $p$ (resp. $P$) without a $+$). Finally, match every yet unmatched string
$2,3,\ldots,n,n-2,\ldots,h'$ in $\bbr$ corresponding to a $-$ in $p$ with a yet unmatched
$1,2,\ldots,h$ corresponding to a column without $+$ in $P$. Again $e_1$ cannot act on
such a combination. This is equivalent to step (3).

Note that $e_1$ acts on the leftmost possible unmatched letter in a tableau and all
$f_i$ act rightmost. Hence $e_1$ will want to transform the smallest unmatched string
$1,2,\ldots,h$ corresponding to a column in $P$ that does not contain a $+$ to a column in $p$. 
This corresponds to the case when an unpaired $+$ exists in $p$ which $e_1$ will move to $P$. 
Else $e_1$ will move the leftmost unpaired $-$ in $P$ to $p$.
\end{proof}

\begin{example}
Take $n>7$ and
\begin{equation*}
P=\young(\blueminus,\mb,\mb\redplus,\mb\mb,\mb\mb\redplus-,\mb\mb\mb\mb)
\quad \text{and} \quad
p=\young(+,\mb,\mb\blueminus,\mb\redplus,\mb\mb\mb\redplus) \; .
\end{equation*}
The paired + of step (1) are circled, as are the paired - of step (2).
There are no pairings for step (3). Hence $e_1(P,p)=(\tilde{P},\tilde{p})$ yields the pair
\begin{equation*}
\tilde{P}=\young(\blueminus,+,\mb\redplus,\mb\mb,\mb\mb\redplus-,\mb\mb\mb\mb) 
\quad \text{and} \quad
\tilde{p}=\young(\mb,\mb\blueminus,\mb\redplus,\mb\mb\mb\redplus) \; .
\end{equation*}
\end{example}

\begin{example}
Take $n>7$ and
\begin{equation*}
P=\young(\blueminus,\mb,\mb\redplus,\mb\mb,\mb\mb\mb+\blueminus,\mb\mb\mb\mb\mb)
\quad \text{and} \quad
p=\young(\greenplus,\mb,\mb\blueminus,\mb\redplus+,\mb\mb\mb\greenminus\blueminus) \; .
\end{equation*}
The paired + of step (1) are circled, as are the paired - of step (2).
The pairings for step (3) are denoted by a square frame. Hence $e_1(P,p)=(\tilde{P},\tilde{p})$ 
yields the pair
\begin{equation*}
P=\young(\blueminus,\mb,\mb\redplus,\mb\mb,\mb\mb++\blueminus,\mb\mb\mb\mb\mb)
\quad \text{and} \quad
p=\young(\greenplus,\mb,\mb\blueminus,\mb\redplus,\mb\mb\mb\greenminus\blueminus) \; .
\end{equation*}
\end{example}

\begin{lemma}
The commutative diagram~\eqref{eq:e1 sigma} holds for $X_{n-2}$ highest weight elements.
\end{lemma}

\begin{proof}
Since $X_{n-2}$ highest weight elements are in bijection with tuples $(P,p)$ of 
$\pm$ diagrams, it suffices to prove~\eqref{eq:e1 sigma} for tuples $(P,p)$.
We need to distinguish the two cases when a $+$ moves from $p$ to $P$ as in~\eqref{eq:+1}
and~\eqref{eq:+2} or when a $-$ moves from $P$ to $p$ as in~\eqref{eq:-1} and~\eqref{eq:-2}.
The following diagrams schematically indicate the region of the $\pm$ diagrams on which $e_1$ and
subsequently $\sigma$ acts. The usual $\pm$ symbols belong to $P$, whereas the circled $\pm$
symbols belong to $p$.
\begin{equation} \label{eq:+1}
\xymatrix{
\scalebox{.8}{
\psset{unit=0.7cm}
\begin{pspicture}(5,6)
\psframe(0,4)(1,5) \rput(0.5,4.5){$-$}
\psframe(0,3)(1,4) \rput(0.5,3.5){$+$}
\psframe(0,2)(1,3) \rput(0.5,2.5){$\redplus$}
\psframe(1,2)(2,3) \rput(1.5,2.5){$\redminus$}
\psframe(2,2)(3,3) \rput(2.5,2.5){$\redminus$}
\psframe(3,2)(4,3) \rput(3.5,2.5){$+$}
\psframe(4,2)(5,3) \rput(4.5,2.5){$-$}
\psframe(2,1)(3,2) \rput(2.5,1.5){$\redplus$}
\psframe(3,1)(4,2) \rput(3.5,1.5){$\redplus$}
\psframe(4,1)(5,2) \rput(4.5,1.5){$\redminus$}
{\red \rput(0.5,1.5){$\underbrace{}_\alpha$}
\rput(1.5,1.5){$\underbrace{}_\beta$}
\rput(2.5,0.5){$\underbrace{}_\gamma$}
\rput(3.5,0.5){$\underbrace{}_\delta$}
\rput(4.5,0.5){$\underbrace{}_\epsilon$}}
\rput(0.5,5.5){$\overbrace{}^a$}
\rput(2,3.5){$\overbrace{\hspace{1.3cm}}^b$}
\rput(3.5,3.5){$\overbrace{}^c$}
\rput(4.5,3.5){$\overbrace{}^d$}
\end{pspicture}}
\qquad
\ar[d]_{\text{$\sigma \circ e_1$ if $\gamma>0$}}
\ar[dr]^{\text{$\sigma \circ e_1$ if $\gamma=0$}} \ar[r]^{\sigma} & 
\qquad
\scalebox{.8}{
\psset{unit=0.7cm}
\begin{pspicture}(5,6)
\psframe(0,4)(1,5) \rput(0.5,4.5){$-$}
\psframe(0,3)(1,4) \rput(0.5,3.5){$+$}
\psframe(0,2)(1,3) \rput(0.5,2.5){$\redplus$}
\psframe(1,2)(2,3) \rput(1.5,2.5){$\redminus$}
\psframe(2,2)(3,3) \rput(2.5,2.5){$\redminus$}
\psframe(3,2)(4,3) \rput(3.5,2.5){$+$}
\psframe(4,2)(5,3) \rput(4.5,2.5){$-$}
\psframe(2,1)(3,2) \rput(2.5,1.5){$\redplus$}
\psframe(3,1)(4,2) \rput(3.5,1.5){$\redplus$}
\psframe(4,1)(5,2) \rput(4.5,1.5){$\redminus$}
{\red \rput(0.5,1.5){$\underbrace{}_\alpha$}
\rput(1.5,1.5){$\underbrace{}_\beta$}
\rput(2.5,0.5){$\underbrace{}_\gamma$}
\rput(3.5,0.5){$\underbrace{}_\delta$}
\rput(4.5,0.5){$\underbrace{}_\epsilon$}}
\rput(0.5,5.5){$\overbrace{}^b$}
\rput(2,3.5){$\overbrace{\hspace{1.3cm}}^a$}
\rput(3.5,3.5){$\overbrace{}^d$}
\rput(4.5,3.5){$\overbrace{}^c$}
\end{pspicture}} \\
\scalebox{.8}{
\psset{unit=0.7cm}
\begin{pspicture}(5,6)
\psframe(0,4)(1,5) \rput(0.5,4.5){$-$}
\psframe(0,3)(1,4) \rput(0.5,3.5){$+$}
\psframe(0,2)(1,3) \rput(0.5,2.5){$\redplus$}
\psframe(1,2)(2,3) \rput(1.5,2.5){$\redminus$}
\psframe(2,2)(3,3) \rput(2.5,2.5){$\redminus$}
\psframe(3,2)(4,3) \rput(3.5,2.5){$+$}
\psframe(4,2)(5,3) \rput(4.5,2.5){$-$}
\psframe(2,1)(3,2) \rput(2.5,1.5){$\redplus$}
\psframe(3,1)(4,2) \rput(3.5,1.5){$\redplus$}
\psframe(4,1)(5,2) \rput(4.5,1.5){$\redminus$}
{\red \rput(0.5,1.5){$\underbrace{}_\alpha$}
\rput(1.5,1.5){$\underbrace{}_\beta$}
\rput(2.5,0.5){$\underbrace{}_{\gamma-1}$}
\rput(3.5,0.5){$\underbrace{}_\delta$}
\rput(4.5,0.5){$\underbrace{}_{\epsilon+1}$}}
\rput(0.5,5.5){$\overbrace{}^{b-1}$}
\rput(2,3.5){$\overbrace{\hspace{1.3cm}}^a$}
\rput(3.5,3.5){$\overbrace{}^d$}
\rput(4.5,3.5){$\overbrace{}^{c+1}$}
\end{pspicture}}
&
\scalebox{.8}{
\psset{unit=0.7cm}
\begin{pspicture}(4,6)
\psframe(0,4)(1,5) \rput(0.5,4.5){$-$}
\psframe(0,3)(1,4) \rput(0.5,3.5){$+$}
\psframe(0,2)(1,3) \rput(0.5,2.5){$\redplus$}
\psframe(1,2)(2,3) \rput(1.5,2.5){$\redminus$}
\psframe(2,2)(3,3) \rput(2.5,2.5){$+$}
\psframe(3,2)(4,3) \rput(3.5,2.5){$-$}
\psframe(2,1)(3,2) \rput(2.5,1.5){$\redplus$}
\psframe(3,1)(4,2) \rput(3.5,1.5){$\redminus$}
{\red \rput(0.5,1.5){$\underbrace{}_{\alpha-1}$}
\rput(1.5,1.5){$\underbrace{}_\beta$}
\rput(2.5,0.5){$\underbrace{}_\delta$}
\rput(3.5,0.5){$\underbrace{}_{\epsilon}$}}
\rput(0.5,5.5){$\overbrace{}^{b-1}$}
\rput(1.5,3.5){$\overbrace{}^a$}
\rput(2.5,3.5){$\overbrace{}^d$}
\rput(3.5,3.5){$\overbrace{}^{c+1}$}
\end{pspicture}}
}
\end{equation}

\begin{equation} \label{eq:+2}
\xymatrix{
\scalebox{.8}{
\psset{unit=0.7cm}
\begin{pspicture}(5,5)
\psframe(0,3)(1,4) \rput(0.5,3.5){$+$}
\psframe(1,3)(2,4) \rput(1.5,3.5){$-$}
\psframe(2,3)(3,4) \rput(2.5,3.5){$-$}
\psframe(3,3)(4,4) \rput(3.5,3.5){$-$}
\psframe(0,2)(1,3) \rput(0.5,2.5){$\redplus$}
\psframe(1,2)(2,3) \rput(1.5,2.5){$\redminus$}
\psframe(2,2)(3,3) \rput(2.5,2.5){$\redminus$}
\psframe(3,2)(4,3) \rput(3.5,2.5){$+$}
\psframe(2,1)(3,2) \rput(2.5,1.5){$\redplus$}
\psframe(3,1)(4,2) \rput(3.5,1.5){$\redplus$}
\psframe(4,1)(5,2) \rput(4.5,1.5){$\redminus$}
{\red \rput(0.5,1.5){$\underbrace{}_\alpha$}
\rput(1.5,1.5){$\underbrace{}_\beta$}
\rput(2.5,0.5){$\underbrace{}_\gamma$}
\rput(3.5,0.5){$\underbrace{}_\delta$}
\rput(4.5,0.5){$\underbrace{}_\epsilon$}}
\rput(0.5,4.5){$\overbrace{}^a$}
\rput(2,4.5){$\overbrace{\hspace{1.3cm}}^b$}
\rput(3.5,4.5){$\overbrace{}^c$}
\rput(4.5,2.5){$\overbrace{}^d$}
\end{pspicture}}
\qquad
\ar[d]_{\text{$\sigma \circ e_1$ if $\gamma>0$}}
\ar[dr]^{\text{$\sigma \circ e_1$ if $\gamma=0$}} \ar[r]^{\sigma} & 
\qquad
\scalebox{.8}{
\psset{unit=0.7cm}
\begin{pspicture}(5,5)
\psframe(0,3)(1,4) \rput(0.5,3.5){$+$}
\psframe(1,3)(2,4) \rput(1.5,3.5){$-$}
\psframe(2,3)(3,4) \rput(2.5,3.5){$-$}
\psframe(3,3)(4,4) \rput(3.5,3.5){$-$}
\psframe(0,2)(1,3) \rput(0.5,2.5){$\redplus$}
\psframe(1,2)(2,3) \rput(1.5,2.5){$\redminus$}
\psframe(2,2)(3,3) \rput(2.5,2.5){$\redminus$}
\psframe(3,2)(4,3) \rput(3.5,2.5){$+$}
\psframe(2,1)(3,2) \rput(2.5,1.5){$\redplus$}
\psframe(3,1)(4,2) \rput(3.5,1.5){$\redplus$}
\psframe(4,1)(5,2) \rput(4.5,1.5){$\redminus$}
{\red \rput(0.5,1.5){$\underbrace{}_\alpha$}
\rput(1.5,1.5){$\underbrace{}_\beta$}
\rput(2.5,0.5){$\underbrace{}_\gamma$}
\rput(3.5,0.5){$\underbrace{}_\delta$}
\rput(4.5,0.5){$\underbrace{}_\epsilon$}}
\rput(0.5,4.5){$\overbrace{}^b$}
\rput(2,4.5){$\overbrace{\hspace{1.3cm}}^a$}
\rput(3.5,4.5){$\overbrace{}^d$}
\rput(4.5,2.5){$\overbrace{}^c$}
\end{pspicture}} \\
\scalebox{.8}{
\psset{unit=0.7cm}
\begin{pspicture}(5,5)
\psframe(0,3)(1,4) \rput(0.5,3.5){$+$}
\psframe(1,3)(2,4) \rput(1.5,3.5){$-$}
\psframe(2,3)(3,4) \rput(2.5,3.5){$-$}
\psframe(3,3)(4,4) \rput(3.5,3.5){$-$}
\psframe(0,2)(1,3) \rput(0.5,2.5){$\redplus$}
\psframe(1,2)(2,3) \rput(1.5,2.5){$\redminus$}
\psframe(2,2)(3,3) \rput(2.5,2.5){$\redminus$}
\psframe(3,2)(4,3) \rput(3.5,2.5){$+$}
\psframe(2,1)(3,2) \rput(2.5,1.5){$\redplus$}
\psframe(3,1)(4,2) \rput(3.5,1.5){$\redplus$}
\psframe(4,1)(5,2) \rput(4.5,1.5){$\redminus$}
{\red \rput(0.5,1.5){$\underbrace{}_\alpha$}
\rput(1.5,1.5){$\underbrace{}_\beta$}
\rput(2.5,0.5){$\underbrace{}_{\gamma-1}$}
\rput(3.5,0.5){$\underbrace{}_\delta$}
\rput(4.5,0.5){$\underbrace{}_{\epsilon+1}$}}
\rput(0.5,4.5){$\overbrace{}^{b-1}$}
\rput(2,4.5){$\overbrace{\hspace{1.3cm}}^a$}
\rput(3.5,4.5){$\overbrace{}^d$}
\rput(4.5,2.5){$\overbrace{}^{c+1}$}
\end{pspicture}}
&
\scalebox{.8}{
\psset{unit=0.7cm}
\begin{pspicture}(4,5)
\psframe(0,3)(1,4) \rput(0.5,3.5){$+$}
\psframe(1,3)(2,4) \rput(1.5,3.5){$-$}
\psframe(2,3)(3,4) \rput(2.5,3.5){$-$}
\psframe(0,2)(1,3) \rput(0.5,2.5){$\redplus$}
\psframe(1,2)(2,3) \rput(1.5,2.5){$\redminus$}
\psframe(2,2)(3,3) \rput(2.5,2.5){$+$}
\psframe(2,1)(3,2) \rput(2.5,1.5){$\redplus$}
\psframe(3,1)(4,2) \rput(3.5,1.5){$\redminus$}
{\red \rput(0.5,1.5){$\underbrace{}_{\alpha-1}$}
\rput(1.5,1.5){$\underbrace{}_\beta$}
\rput(2.5,0.5){$\underbrace{}_\delta$}
\rput(3.5,0.5){$\underbrace{}_{\epsilon}$}}
\rput(0.5,4.5){$\overbrace{}^{b-1}$}
\rput(1.5,4.5){$\overbrace{}^a$}
\rput(2.5,4.5){$\overbrace{}^d$}
\rput(3.5,2.5){$\overbrace{}^{c+1}$}
\end{pspicture}}
}
\end{equation}

\begin{equation} \label{eq:-1}
\xymatrix{
\scalebox{.8}{
\psset{unit=0.7cm}
\begin{pspicture}(5,5)
\psframe(0,3)(1,4) \rput(0.5,3.5){$+$}
\psframe(1,3)(2,4) \rput(1.5,3.5){$-$}
\psframe(2,3)(3,4) \rput(2.5,3.5){$-$}
\psframe(3,3)(4,4) \rput(3.5,3.5){$-$}
\psframe(0,2)(1,3) \rput(0.5,2.5){$\redplus$}
\psframe(1,2)(2,3) \rput(1.5,2.5){$\redminus$}
\psframe(2,2)(3,3) \rput(2.5,2.5){$\redminus$}
\psframe(3,2)(4,3) \rput(3.5,2.5){$+$}
\psframe(2,1)(3,2) \rput(2.5,1.5){$\redplus$}
\psframe(3,1)(4,2) \rput(3.5,1.5){$\redplus$}
\psframe(4,1)(5,2) \rput(4.5,1.5){$\redminus$}
{\red \rput(0.5,1.5){$\underbrace{}_\alpha$}
\rput(1.5,1.5){$\underbrace{}_\beta$}
\rput(2.5,0.5){$\underbrace{}_\gamma$}
\rput(3.5,0.5){$\underbrace{}_\delta$}
\rput(4.5,0.5){$\underbrace{}_\epsilon$}}
\rput(0.5,4.5){$\overbrace{}^a$}
\rput(2,4.5){$\overbrace{\hspace{1.3cm}}^b$}
\rput(3.5,4.5){$\overbrace{}^c$}
\rput(4.5,2.5){$\overbrace{}^d$}
\end{pspicture}}
\qquad
\ar[d]_{\text{$\sigma \circ e_1$ if $\delta+\epsilon=c+d$}}
\ar[dr]^{\text{$\sigma \circ e_1$ if $\delta+\epsilon<c+d$}} \ar[r]^{\sigma} & 
\qquad
\scalebox{.8}{
\psset{unit=0.7cm}
\begin{pspicture}(5,5)
\psframe(0,3)(1,4) \rput(0.5,3.5){$+$}
\psframe(1,3)(2,4) \rput(1.5,3.5){$-$}
\psframe(2,3)(3,4) \rput(2.5,3.5){$-$}
\psframe(3,3)(4,4) \rput(3.5,3.5){$-$}
\psframe(0,2)(1,3) \rput(0.5,2.5){$\redplus$}
\psframe(1,2)(2,3) \rput(1.5,2.5){$\redminus$}
\psframe(2,2)(3,3) \rput(2.5,2.5){$\redminus$}
\psframe(3,2)(4,3) \rput(3.5,2.5){$+$}
\psframe(2,1)(3,2) \rput(2.5,1.5){$\redplus$}
\psframe(3,1)(4,2) \rput(3.5,1.5){$\redplus$}
\psframe(4,1)(5,2) \rput(4.5,1.5){$\redminus$}
{\red \rput(0.5,1.5){$\underbrace{}_\alpha$}
\rput(1.5,1.5){$\underbrace{}_\beta$}
\rput(2.5,0.5){$\underbrace{}_\gamma$}
\rput(3.5,0.5){$\underbrace{}_\delta$}
\rput(4.5,0.5){$\underbrace{}_\epsilon$}}
\rput(0.5,4.5){$\overbrace{}^b$}
\rput(2,4.5){$\overbrace{\hspace{1.3cm}}^a$}
\rput(3.5,4.5){$\overbrace{}^d$}
\rput(4.5,2.5){$\overbrace{}^c$}
\end{pspicture}} \\
\scalebox{.8}{
\psset{unit=0.7cm}
\begin{pspicture}(5,5)
\psframe(0,3)(1,4) \rput(0.5,3.5){$+$}
\psframe(1,3)(2,4) \rput(1.5,3.5){$-$}
\psframe(2,3)(3,4) \rput(2.5,3.5){$-$}
\psframe(3,3)(4,4) \rput(3.5,3.5){$-$}
\psframe(0,2)(1,3) \rput(0.5,2.5){$\redplus$}
\psframe(1,2)(2,3) \rput(1.5,2.5){$\redminus$}
\psframe(2,2)(3,3) \rput(2.5,2.5){$\redminus$}
\psframe(3,2)(4,3) \rput(3.5,2.5){$+$}
\psframe(2,1)(3,2) \rput(2.5,1.5){$\redplus$}
\psframe(3,1)(4,2) \rput(3.5,1.5){$\redplus$}
\psframe(4,1)(5,2) \rput(4.5,1.5){$\redminus$}
{\red \rput(0.5,1.5){$\underbrace{}_\alpha$}
\rput(1.5,1.5){$\underbrace{}_\beta$}
\rput(2.5,0.5){$\underbrace{}_{\gamma+1}$}
\rput(3.5,0.5){$\underbrace{}_{\delta-1}$}
\rput(4.5,0.5){$\underbrace{}_\epsilon$}}
\rput(0.5,4.5){$\overbrace{}^{b}$}
\rput(2,4.5){$\overbrace{\hspace{1.3cm}}^{a+1}$}
\rput(3.5,4.5){$\overbrace{}^d$}
\rput(4.5,2.5){$\overbrace{}^{c-1}$}
\end{pspicture}}
& 
\scalebox{.8}{
\psset{unit=0.7cm}
\begin{pspicture}(4,5)
\psframe(0,3)(1,4) \rput(0.5,3.5){$+$}
\psframe(1,3)(2,4) \rput(1.5,3.5){$-$}
\psframe(2,3)(3,4) \rput(2.5,3.5){$-$}
\psframe(0,2)(1,3) \rput(0.5,2.5){$\redplus$}
\psframe(1,2)(2,3) \rput(1.5,2.5){$\redminus$}
\psframe(2,2)(3,3) \rput(2.5,2.5){$+$}
\psframe(2,1)(3,2) \rput(2.5,1.5){$\redplus$}
\psframe(3,1)(4,2) \rput(3.5,1.5){$\redminus$}
{\red \rput(0.5,1.5){$\underbrace{}_\alpha$}
\rput(1.5,1.5){$\underbrace{}_{\beta+1}$}
\rput(2.5,0.5){$\underbrace{}_\delta$}
\rput(3.5,0.5){$\underbrace{}_\epsilon$}}
\rput(0.5,4.5){$\overbrace{}^{b}$}
\rput(1.5,4.5){$\overbrace{}^{a+1}$}
\rput(2.5,4.5){$\overbrace{}^d$}
\rput(3.5,2.5){$\overbrace{}^{c-1}$}
\end{pspicture}}
}
\end{equation}

\begin{equation} \label{eq:-2}
\xymatrix{
\scalebox{.8}{
\psset{unit=0.7cm}
\begin{pspicture}(5,6)
\psframe(0,4)(1,5) \rput(0.5,4.5){$-$}
\psframe(0,3)(1,4) \rput(0.5,3.5){$+$}
\psframe(0,2)(1,3) \rput(0.5,2.5){$\redplus$}
\psframe(1,2)(2,3) \rput(1.5,2.5){$\redminus$}
\psframe(2,2)(3,3) \rput(2.5,2.5){$\redminus$}
\psframe(3,2)(4,3) \rput(3.5,2.5){$+$}
\psframe(4,2)(5,3) \rput(4.5,2.5){$-$}
\psframe(2,1)(3,2) \rput(2.5,1.5){$\redplus$}
\psframe(3,1)(4,2) \rput(3.5,1.5){$\redplus$}
\psframe(4,1)(5,2) \rput(4.5,1.5){$\redminus$}
{\red \rput(0.5,1.5){$\underbrace{}_\alpha$}
\rput(1.5,1.5){$\underbrace{}_\beta$}
\rput(2.5,0.5){$\underbrace{}_\gamma$}
\rput(3.5,0.5){$\underbrace{}_\delta$}
\rput(4.5,0.5){$\underbrace{}_\epsilon$}}
\rput(0.5,5.5){$\overbrace{}^a$}
\rput(2,3.5){$\overbrace{\hspace{1.3cm}}^b$}
\rput(3.5,3.5){$\overbrace{}^c$}
\rput(4.5,3.5){$\overbrace{}^d$}
\end{pspicture}}
\qquad
\ar[d]_{\text{$\sigma \circ e_1$ if $\delta+\epsilon=c+d$}}
\ar[dr]^{\text{$\sigma \circ e_1$ if $\delta+\epsilon<c+d$}} \ar[r]^{\sigma} & 
\qquad
\scalebox{.8}{
\psset{unit=0.7cm}
\begin{pspicture}(5,6)
\psframe(0,4)(1,5) \rput(0.5,4.5){$-$}
\psframe(0,3)(1,4) \rput(0.5,3.5){$+$}
\psframe(0,2)(1,3) \rput(0.5,2.5){$\redplus$}
\psframe(1,2)(2,3) \rput(1.5,2.5){$\redminus$}
\psframe(2,2)(3,3) \rput(2.5,2.5){$\redminus$}
\psframe(3,2)(4,3) \rput(3.5,2.5){$+$}
\psframe(4,2)(5,3) \rput(4.5,2.5){$-$}
\psframe(2,1)(3,2) \rput(2.5,1.5){$\redplus$}
\psframe(3,1)(4,2) \rput(3.5,1.5){$\redplus$}
\psframe(4,1)(5,2) \rput(4.5,1.5){$\redminus$}
{\red \rput(0.5,1.5){$\underbrace{}_\alpha$}
\rput(1.5,1.5){$\underbrace{}_\beta$}
\rput(2.5,0.5){$\underbrace{}_\gamma$}
\rput(3.5,0.5){$\underbrace{}_\delta$}
\rput(4.5,0.5){$\underbrace{}_\epsilon$}}
\rput(0.5,5.5){$\overbrace{}^b$}
\rput(2,3.5){$\overbrace{\hspace{1.3cm}}^a$}
\rput(3.5,3.5){$\overbrace{}^d$}
\rput(4.5,3.5){$\overbrace{}^c$}
\end{pspicture}} \\
\scalebox{.8}{
\psset{unit=0.7cm}
\begin{pspicture}(5,6)
\psframe(0,4)(1,5) \rput(0.5,4.5){$-$}
\psframe(0,3)(1,4) \rput(0.5,3.5){$+$}
\psframe(0,2)(1,3) \rput(0.5,2.5){$\redplus$}
\psframe(1,2)(2,3) \rput(1.5,2.5){$\redminus$}
\psframe(2,2)(3,3) \rput(2.5,2.5){$\redminus$}
\psframe(3,2)(4,3) \rput(3.5,2.5){$+$}
\psframe(4,2)(5,3) \rput(4.5,2.5){$-$}
\psframe(2,1)(3,2) \rput(2.5,1.5){$\redplus$}
\psframe(3,1)(4,2) \rput(3.5,1.5){$\redplus$}
\psframe(4,1)(5,2) \rput(4.5,1.5){$\redminus$}
{\red \rput(0.5,1.5){$\underbrace{}_\alpha$}
\rput(1.5,1.5){$\underbrace{}_\beta$}
\rput(2.5,0.5){$\underbrace{}_{\gamma+1}$}
\rput(3.5,0.5){$\underbrace{}_{\delta-1}$}
\rput(4.5,0.5){$\underbrace{}_\epsilon$}}
\rput(0.5,5.5){$\overbrace{}^{b+1}$}
\rput(2,3.5){$\overbrace{\hspace{1.3cm}}^a$}
\rput(3.5,3.5){$\overbrace{}^{d-1}$}
\rput(4.5,3.5){$\overbrace{}^c$}
\end{pspicture}}
&
\scalebox{.8}{
\psset{unit=0.7cm}
\begin{pspicture}(4,6)
\psframe(0,4)(1,5) \rput(0.5,4.5){$-$}
\psframe(0,3)(1,4) \rput(0.5,3.5){$+$}
\psframe(0,2)(1,3) \rput(0.5,2.5){$\redplus$}
\psframe(1,2)(2,3) \rput(1.5,2.5){$\redminus$}
\psframe(2,2)(3,3) \rput(2.5,2.5){$+$}
\psframe(3,2)(4,3) \rput(3.5,2.5){$-$}
\psframe(2,1)(3,2) \rput(2.5,1.5){$\redplus$}
\psframe(3,1)(4,2) \rput(3.5,1.5){$\redminus$}
{\red \rput(0.5,1.5){$\underbrace{}_\alpha$}
\rput(1.5,1.5){$\underbrace{}_{\beta+1}$}
\rput(2.5,0.5){$\underbrace{}_{\delta}$}
\rput(3.5,0.5){$\underbrace{}_{\epsilon}$}}
\rput(0.5,5.5){$\overbrace{}^{b+1}$}
\rput(1.5,3.5){$\overbrace{}^a$}
\rput(2.5,3.5){$\overbrace{}^{d-1}$}
\rput(3.5,3.5){$\overbrace{}^c$}
\end{pspicture}}
}
\end{equation}
Let $(P,p)$ be the $\pm$ diagram on the left top in each case. We show that
$e_1$ acts in the same way on both the top right diagram $\sigma(P,p)$ and the bottom 
diagram $\sigma\circ e_1(P,p)$, except for a few special cases which are treated 
separately below. Up to the special cases, this proves~\eqref{eq:e1 sigma}. 
For case~\eqref{eq:+1} with $\gamma>0$, the net pairing according to the algorithm
before Lemma~\ref{lem:e1 action} is the same for the top right and bottom diagram, 
except in the case when $e_1$ acts on one of the $\gamma$ pluses in $\sigma(P,p)$ 
and $\gamma=1$. Let us call this case S:\eqref{eq:+1} $\gamma=1$ which will be 
treated separately later.
For~\eqref{eq:+1} with $\gamma=0$, the condition that a plus moves
from $p$ to $P$ requires that $\alpha>a$. Since $\alpha+\beta\le a+b$ we have 
$\beta<b$. Hence the only special case that needs to be checked explicitly, occurs when 
$e_1$ acts on the extra minus in the set of $c+1$ minuses in $\sigma \circ e_1(P,p)$.
Call this case S:\eqref{eq:+1} $\gamma=0$.
For case~\eqref{eq:+2} with $\gamma=0$ the pairing is obviously the same for $\sigma(P,p)$
and $\sigma\circ e_1(P,p)$. For $\gamma>0$, the only special case occurs when $e_1$ acts 
on one of the $\gamma$ pluses in $\sigma(P,p)$ and $\gamma=1$. Call this case
S:\eqref{eq:+2} $\gamma=1$.

In case~\eqref{eq:-1} with $\delta+\epsilon=c+d$, the fact that a minus moves from $P$ to
$p$ requires that $c>\epsilon$ so that $\delta>d$. Hence the action of $e_1$ is the same on
$\sigma(P,p)$ and $\sigma\circ e_1(P,p)$ except for the special case when $\delta=d+1$
and $e_1$ acts on one of the $\gamma+1$ pluses in $\sigma\circ e_1(P,p)$ and 
on one of the $\delta$ pluses in $\sigma(P,p)$. Call this case S:\eqref{eq:-1} $\delta=d+1$.
For~\eqref{eq:-1} with $\delta+\epsilon<c+d$, the action of $e_1$ on $\sigma\circ e_1(P,p)$ and 
$\sigma(P,p)$ is the same.
Similarly, for~\eqref{eq:-2}, a minus in the group of $d$ minuses in $P$ must move under
the action of $e_1$. This requires $\beta\ge a$ or equivalently, since $\alpha+\beta\le a+b$, it 
follows that $\alpha\le b$. Similarly, we need $\epsilon<d$ and $\delta>c$. Hence the net 
pairing remains unchanged and $e_1$ acts in the same way on $\sigma\circ e_1(P,p)$ and 
$e_1(P,p)$, except possibly when $e_1$ acts on one of the $\delta$ pluses in $\sigma(P,p)$, 
but on one of the $\gamma+1$ pluses in $\sigma\circ e_1(P,p)$. But for $\delta+\epsilon=c+d$
this is impossible.

We are now going to treat the special cases. For S:\eqref{eq:+1} $\gamma=1$, 
case~\eqref{eq:+1} $\gamma=1$ applies to $\sigma(P,p)$, and case~\eqref{eq:+1} 
$\gamma=0$ applies to both $\sigma\circ e_1(P,p)$ and $\sigma\circ e_1\circ\sigma(P,p)$.
It can be shown explicitly that~\eqref{eq:e1 sigma} commutes. The same happens for
S:\eqref{eq:+2} $\gamma=1$ with \eqref{eq:+1} replaced by \eqref{eq:+2} everywhere.
For S:\eqref{eq:+1} $\gamma=0$, $e_1$ acts on one of the $c+1$ minuses in $\sigma\circ e_1(P,p)$.
Then case~\eqref{eq:-2} with $\delta>0$ applies for $\sigma\circ e_1(P,p)$, case~\eqref{eq:-1} 
with $\delta>0$ applies for $\sigma(P,p)$ and case~\eqref{eq:+2} with $\gamma>0$ applies for 
$\sigma\circ e_1 \circ \sigma(P,p)$. Again it can be shown explicitly that~\eqref{eq:e1 sigma} 
commutes.

For S:\eqref{eq:-1} $\delta=d+1$, $e_1$ acts on the group of $\gamma+1$ pluses in
$\sigma\circ e_1(P,p)$, but the group of $\delta$ pluses in $\sigma(P,p)$. Then 
case~\eqref{eq:+2} with $\gamma>0$ applies for $\sigma\circ e_1(P,p)$,
case~\eqref{eq:+1} with $\gamma=0$ applies for $\sigma(P,p)$, and
case~\eqref{eq:+2} with $\delta+\epsilon<c+d$ applies for $\sigma\circ e_1\circ \sigma(P,p)$.
Again it follows that~\eqref{eq:e1 sigma} commutes.
\end{proof}

\section{Perfectness} \label{sec:perfect}

In this section we show that the KR crystals $B^{r,s}$ of type $D_n^{(1)}$, $B_n^{(1)}$, and 
$A_{2n-1}^{(2)}$ are perfect.
The notion of perfect crystal was introduced in~\cite{KKMMNN:1992} in order to construct
a path realization of highest weight modules. In a subsequent paper~\cite{KKMMNN:1992a}
a perfect crystal of arbitrary level was given for every nonexceptional affine algebra.
Further perfect crystals were found and studied in~\cite{BFKL:2006, HN:2006, JMO:2000,
Ka:2002, Ko:1999, KMOY:2006, NS:2005, Y:1998} for example.
In~\cite{HKOTT:2002,HKOTY:1999} it is conjectured that $B^{r,s}$ is perfect of level $s/t_r$ 
if $s/t_r$ is an integer, where $t_r=\max(1,2/(\alpha_r,\alpha_r))$. 

To define perfect crystals, we need a few preliminary definitions.
Define the set of level $\ell$ weights to be
$(P_{\cl}^{+})_{\ell}=\{\la\in P_{\cl}^{+} \mid \langle c,\la \rangle=\ell\}$.  
For a crystal basis $B$, we define $B_{\min}$ to be the set of crystal basis
elements $b$ such that $\langle c,\ve(b)\rangle$ is minimal over $b\in B$.
We say that $B\otimes B$ is connected if it contains only one irreducible component.

\begin{definition} \label{def:perfect}
A crystal $B$ is a perfect crystal of level $\ell$ if:
\begin{enumerate}
    \item $B\otimes B$ is connected;
    \item there exists $\la\in P_{\cl}$ such that $\wt(B)\subset
    \la+\sum_{i\neq0}\mathbb{Z}_{\leq0}\alpha_{i}$ and
    $\#(B_{\la})=1$;
    \item there is a finite-dimensional irreducible $U_q(\geh)$-module $V$ with a
    crystal base whose crystal graph is isomorphic to $B$;
    \item for any $b\in B$, we have $\langle c,\ve(b)\rangle \geq
    \ell$;
    \item the maps $\ve$ and $\varphi$ from $B_{\min}$ to
    $(P_{\cl}^{+})_{\ell}$ are bijective.
\end{enumerate}
We use the notation $\lev(B)$ to indicate the level of the perfect 
crystal $B$.
\end{definition}

We show here that $B^{r,s}$ of type $D_n^{(1)}$, $B_n^{(1)}$, and $A_{2n-1}^{(2)}$ satisfies conditions
(1), (2), (4), and (5) of Definition~\ref{def:perfect} with level $\ell=s$. Condition (3) is proved 
in~\cite{OS:2007}.
Let us first give an explicit construction of all elements in $B^{r,s}_{\min}$.
Let $\Lambda\in (P_\cl^+)_s$ be a dominant weight of level $s$. That is 
$\Lambda=\ell_0\Lambda_0+\ell_1\Lambda_1+\cdots+\ell_n \Lambda_n$ with 
$\lev(\La):=\langle c, \La \rangle = s$, where explicitly
\begin{align*}
	\lev(\La) &= \ell_0+\ell_1+2\ell_2+2\ell_3+\cdots+2\ell_{n-2}+\ell_{n-1}+\ell_n && \text{for type $D_n^{(1)}$}\\
	\lev(\La) &= \ell_0+\ell_1+2\ell_2+2\ell_3+\cdots+2\ell_{n-2}+2\ell_{n-1}+\ell_n && \text{for type $B_n^{(1)}$}\\
	\lev(\La) &= \ell_0+\ell_1+2\ell_2+2\ell_3+\cdots+2\ell_{n-2}+2\ell_{n-1}+2\ell_n && \text{for type $A_{2n-1}^{(2)}$.}
\end{align*}	 

To a given fundamental weight $\Lambda_k$ we may associate the following $\pm$ diagram
\begin{equation} \label{eq:La pm}
\diagram: \Lambda_k \mapsto \begin{cases}
	\phantom{k+1\left\{\right.} \emptyset & \text{if $r$ is even and $k=0$}\\[1mm]
	\phantom{k+1\left\{\right.} \young(-,+) & \text{if $r$ is even and $k=1$}\\[3mm]
	\phantom{k+1\left\{\right.} \young(+) & \text{if $r$ is odd and $k=0$}\\[2mm]
	\phantom{k+1\left\{\right.} \young(-) & \text{if $r$ is odd and $k=1$}\\[2mm]
	k+1\left\{\young(-,+,\mb\mb,\mb\mb)\right. 
		& \text{if $k\not \equiv r \bmod{2}$ and $2\le k\le r$}\\[8mm]
	\qquad k\left\{\young(+-,\mb\mb,\mb\mb,\mb\mb)\right.
		& \text{if $k \equiv r \bmod{2}$ and $2\le k\le r$}\\[8mm]
	\qquad r\left\{\young(\mb\mb,\mb\mb,\mb\mb,\mb\mb)\right.
		& \begin{array}{l} \text{if $r<k\le n-2$ for type $D_n^{(1)}$}\\
		                                \text{if $r<k\le n-1$ for type $B_n^{(1)}$}\\
		                                \text{if $r<k\le n$ for type $A_{2n-1}^{(2)}$} \end{array} \\[8mm]
	\qquad r\left\{\young(\mb,\mb,\mb,\mb)\right. 
		& \begin{array}{l} \text{if $k=n-1,n$ for type $D_n^{(1)}$}\\
		                                \text{if $k=n$ for type $B_n^{(1)}$.} \end{array}
\end{cases}
\end{equation}
This map can be extended to any dominant weight $\Lambda=\ell_0\Lambda_0+\cdots+
\ell_n\Lambda_n$ by concatenating the columns of the $\pm$ diagrams of each piece. 
\begin{example}\label{ex:minimal}
Let $s=9$ and $\Lambda=\Lambda_0+2\Lambda_1+\Lambda_2+\Lambda_3+\Lambda_5$ of
type $D_8^{(1)}$. Then 
\begin{equation*}
\begin{split}
\diagram(\Lambda) &= \young(\mb\mb+--,\mb\mb\mb\mb+,\mb\mb\mb\mb\mb\mb+--)
	\quad \text{for $r=3$ and}\\
\diagram(\Lambda) &= \young(\mb\mb-,\mb\mb+,\mb\mb\mb\mb+---,\mb\mb\mb\mb\mb\mb++) 
	\quad \text{for $r=4$.}
\end{split}
\end{equation*}
\end{example}

To every fundamental weight $\Lambda_k$ we also associate a string of operators $f_i$ with 
$i\in\{2,3,\ldots,n\}$ as follows. Let $T(\Lambda_k)$ be the tableau assigned to $\Lambda_k$ as
\begin{equation*}
T(\Lambda_k) = \begin{cases}
	\quad \emptyset & \text{if $r$ is even and $k=0$}\\[1mm]
	\quad \young(\bb,2) & \text{if $r$ is even and $k=1$}\\[4mm]
	\quad \young(1) & \text{if $r$ is odd and $k=0$}\\[1mm]
	\quad \young(\ab) & \text{if $r$ is odd and $k=1$}\\[2mm]
	\quad \begin{array}{|c|c|} \cline{1-1} \overline{k+1} & \multicolumn{1}{c}{} \\ \cline{1-1} 
                                        k+1 & \multicolumn{1}{c}{}\\ \hline k&\bar{2}\\
                                        \hline \vdots & \vdots\\ \hline 2& \bar{k}\\
                                         \hline \end{array} & \text{if $2\le k\le r$ and $k\not \equiv r\bmod 2$} \\[1.4cm]
          \quad \begin{array}{|c|c|} \hline k & \ab\\ \hline \vdots & \vdots\\[2mm] \hline 1& \bar{k}\\
                                         \hline \end{array} & \text{if $2\le k\le r$ and $k\equiv r\bmod 2$} \\[1cm]
	\quad \begin{array}{|c|c|} \hline k&\overline{k-r+1}\\ \hline \vdots&\vdots\\ \hline
                                        k-r+1 & \overline{k} \\ \hline \end{array} 
                     & \begin{array}{l} \text{if $r<k\le n-2$ for type $D_n^{(1)}$}\\
                     			      \text{if $r<k\le n-1$ for type $B_n^{(1)}$}\\
			      		      \text{if $r<k\le n$ for type $A_{2n-1}^{(2)}$} \end{array} \\[0.9cm]
          \quad \left.\begin{array}{|c|} \hline \vdots \\ \hline n\\ \hline 
                                         \overline{n}\\ \hline n \\ \hline\end{array} \right\} r 
                     & \text{for $k=n-1$ for type $D_n^{(1)}$}\\[1.1cm]
          \quad \text{previous case with $n\leftrightarrow \bar{n}$}
                & \text{for $k=n$ for type $D_n^{(1)}$}\\[.2cm]
           \quad \left.\begin{array}{|c|} \hline 0 \\ \hline \vdots \\ \hline 0 \\ \hline \end{array} \right\} r 
                     & \text{for $k=n$ for type $B_n^{(1)}$}
\end{cases}
\end{equation*}
Then $f(\Lambda_k)$ for $0\le k\le n$ is defined such that
$T(\Lambda_k) = f(\Lambda_k) Y(\Lambda_k)$, where $Y(\Lambda_k)$ is the $X_{n-1}$-highest
weight tableau corresponding to $\diagram(\Lambda_k)$. Note that in fact $f(\Lambda_0)
=f(\Lambda_1)=1$.

The minimal element $b$ in $B^{r,s}$ that satisfies $\varepsilon(b)=\Lambda$ can now be
constructed as follows
\begin{equation*}
	b = f(\Lambda_n)^{\ell_n} \cdots  f(\Lambda_2)^{\ell_2} \bij(\diagram(\Lambda)).
\end{equation*}

{}From the condition that $\wt(b)=\vp(b)-\ve(b)$ it is not hard to see that $\vp(b)=\ve(b)$
for $b\in B_{\min}$ and $r$ even. For $r$ odd, we have $\vp(b)=\sigma \sigma' \ve(b)$
for $b\in B_{\min}$, where $\sigma$ is the Dynkin diagram automorphism interchanging nodes
0 and 1 and $\sigma'$ is the Dynkin diagram automorphism interchanging nodes $n-1$ and $n$.

\begin{example}
Continuing Example~\ref{ex:minimal}, we find that the element $b\in B^{r,s}$ such that
$\varepsilon(b)=\Lambda$ is
\begin{equation*}
\begin{split}
	&\young(35\eb\ab\ab,225\eb\bb,1115\eb\cb\bb\ab\ab)  \quad \text{for $r=3$ and}\\
	&\young(5\eb\bb,35\eb,225\eb\bb\ab\ab\ab,11125\eb\cb\bb)  \qquad \text{for $r=4$.}
\end{split}
\end{equation*}
In MuPAD-Combinat these examples can be reproduced via
\begin{verbatim}
>> KR:=combinat::crystals::kirillovReshetikhin(3,9,["D",8,1]):
>> KR::minimalElement([1,2,1,1,0,1,0,0,0])
\end{verbatim}
and 
\begin{verbatim}
>> KR:=combinat::crystals::kirillovReshetikhin(4,9,["D",8,1]):
>> KR::minimalElement([1,2,1,1,0,1,0,0,0])
\end{verbatim}
where $\La$ is entered via its coordinate vector in the basis of $\La_i$.
\end{example}

\begin{example}
The list of all minimal elements of $B^{2,2}$ of type $D_4^{(1)}$ is given by
\begin{equation*}
\emptyset \quad \young(4,\db) \quad \young(\db,4) \quad \young(\bb,2) \quad 
\young(\bb\ab,12) \quad \young(2\ab,1\bb) \quad \young(4\ab,1\db) \quad 
\young(\db\ab,14) \quad \young(3\bb,2\cb) \quad \young(4\cb,3\db) \quad
\young(\db\cb,34) \quad .
\end{equation*}
In MuPAD-Combinat the minimal elements are listed via
\begin{verbatim}
>> KR := combinat::crystals::kirillovReshetikhin(2,2,["D",4,1]):
>> KR::minimal()
\end{verbatim}
\end{example}

\begin{theorem}
The crystals $B^{r,s}$ of type $D_n^{(1)}$, $B_n^{(1)}$, and $A_{2n-1}^{(2)}$ are perfect of level $s$.
\end{theorem}

\begin{proof}
Conditions (1) and (2) of Definition~\ref{def:perfect} follow from \cite{FSS:2006}. Condition (3) 
is proven in~\cite[Theorem 1.2]{OS:2007}. We have shown that the
crystals $B^{r,s}$ of type $D_n^{(1)}$, $B_n^{(1)}$, and $A_{2n-1}^{(2)}$ constructed in this paper 
satisfy~\cite[Assumption 1]{FSS:2006}. Hence by~\cite[Corollary 6.1]{FSS:2006}
the tensor product $B^{r,s}\otimes B^{r,s}$ is connected. Furthermore, it follows 
from~\cite[Theorem 4.4]{FSS:2006} that $\lambda=s\omega_r$ is the weight required 
for property (2).

We prove conditions (4) and (5) by induction on $s$. Suppose the classical
dominant weight $\om$ is such that $B(\om) \subset B^{r,s-1}$ according 
to~\eqref{eq:classical decomp}. Since also $B(\om)\subset B^{r,s}$
this defines an embedding $\iota:B^{r,s-1} \hookrightarrow B^{r,s}$. The classical
component remains unchanged under this embedding, so that $\ve_i(b) = \ve_i(\iota(b))$
for all $b\in B^{r,s-1}$ and $i=1,2,\ldots,n$. Furthermore, $\ve_0(\iota(b)) = \ve_0(b)+1$
for $b\in B^{r,s-1}$ due to the fact that $\sigma$ adds an extra $\young(-,+)$ to the corresponding
$\pm$ diagram. Hence by induction $\langle c, \ve(b) \rangle \ge s$ for $b\in B(\om)\subset B^{r,s}$ 
with width of $\om$ strictly less than $s$. 

For (4) it remains to show that $\langle c, \ve(b) \rangle \ge s$ for $b\in B(\om)\subset B^{r,s}$ 
with $\om$ of width $s$. Suppose that $b$ is in the $X_{n-2}$-highest weight component
corresponding to the pair of $\pm$ diagrams $(P,p)$. If $p$ is not of width $s$, then $P$
has at least one column of the form $\young(-,+)$, $\young(-)$, or $\young(+)$. In the first
case $\sigma(b)\in B(\om)\subset B^{r,s}$ with $\om$ of width strictly less than $s$. Hence 
by the same arguments as above $\langle c, \ve(\sigma(b)) \rangle = \langle c, \ve(b) \rangle \ge s$.
If $P$ contains a column $\young(-)$, then the last column of $b$ is $\young(\ab)$. In
this case $b$ is in the image of the embedding $\tilde{\iota} : B^{r,s-1} \hookrightarrow B^{r,s}$ for
$r$ odd, which adds a column $\young(\ab)$. Under this embedding we have $\ve_i(\tilde{b}) 
= \ve_i(\tilde{\iota}(\tilde{b}))$ for $i =0,2,3,\ldots,n$ and $\ve_1(\iota(\tilde{b})) = \ve_1(\tilde{b})+1$
for all $\tilde{b}\in B^{r,s-1}$. Hence by induction on $s$, we have $\langle c, \ve(b) \rangle \ge s$.
Similarly, if $P$ contains a column $\young(+)$, then $\sigD(P)$ contains a column $\young(-)$ and by 
the above arguments $\langle c, \ve(b) \rangle = \langle c, \ve(\sigma(b)) \rangle \ge s$.
Hence from now on we may assume that $p$ has width $s$.

Let $k_+, k_-, k_\mp$,  and $k_0$ be the number of columns in $p$ containing only a $+, -, \mp$, and 
no $+$ or $-$, respectively. Then we have
\begin{equation*}
	\left.
	\begin{array}{ll} 
		\sum_{i=2}^{n-2} 2\ve_i(b) + \ve_{n-1}(b) + \ve_n(b) & \text{for type $D_n^{(1)}$}\\
		\sum_{i=2}^{n-1} 2\ve_i(b) + \ve_n(b) & \text{for type $B_n^{(1)}$}\\
		\sum_{i=2}^n 2\ve_i(b) & \text{for type $A_{2n-1}^{(2)}$}
	\end{array}
	\right\}
	\ge 	2k_- + k_\mp + k_0 = s + s' - 2k_+
\end{equation*}
where $s'=s-k_0-k_\mp$. This inequality follows from the fact that for an $X_{n-2}$ highest weight 
vector corresponding to $p$ each $-$ and each column without a $+$ contributes 2 to the sum on 
the left-hand side, which is $(2k_- +2 k_{\mp}) + (2 k_0 + 2 k_-)=4k_- + 2 k_\mp + 2 k_0$. 
Using products of $f_i$ it is possible to either change a $\ve_2$ to a $\ve_{n-1}$ or $\ve_n$ 
for type $D_n^{(1)}$ and $B_n^{(1)}$, or to match two contributions in two different columns
(an $\bar{i}$ can support a $i$ in the tableau). In both cases the contribution can at most be halved.

By Lemma~\ref{lem:e1 action} we obtain
\begin{equation*}
\begin{split}
	\ve_1(b) &\ge \max(k_+-K_+ + K_- - k_-,0),\\
	\ve_0(b) &\ge \max(k_+-K_- + K_+ - k_-,0),
\end{split}
\end{equation*}
where $K_+$ and $K_-$ are the total number of $+$ and $-$ in $P$. Hence
\begin{equation} \label{eq:lev ineq}
	\langle c, \ve(b) \rangle \ge 
	s+s'-2k_++2\max(k_+-k_-,0).
\end{equation}
If $k_+\le s'/2$, then certainly $\langle c, \ve(b) \rangle \ge s$ by~\eqref{eq:lev ineq}. 
If $k_+>s'/2$, then $k_-<s'/2$ since $k_+ + k_- = s'$, so that by~\eqref{eq:lev ineq} we have
$\langle c, \ve(b) \rangle \ge s+s'-2k_->s$. This concludes the proof of property (4).

Finally let us prove property (5).
For all weights $\La\in (P_{\cl}^{+})_s$  which contain $\Lambda_0$ (resp. $\Lambda_1$), we know 
by induction on $s$ that $\La-\La_0$ (resp. $\La-\La_1$) is in bijection  with a minimal element 
$b\in B^{r,s-1}$, and hence the weight $\La$ corresponds to $\iota(b)$ or $\tilde{\iota}(\sigma(b))$ 
(resp. $\iota(\sigma(b))$ or $\tilde{\iota}(b)$) which is a minimal element in $B^{r,s}$.
Hence we may restrict our attention to $b\in B^{r,s}$ in an $X_{n-2}$ highest weight component 
labeled by a tuple of $\pm$ diagrams $(P,p)$ such that $p$ has width $s$. By the arguments 
following~\eqref{eq:lev ineq} it is clear that $\langle c, \ve(b) \rangle = s$ is only possible if 
$k_+=s'/2=k_-$. In this case also $K_+=K_-$, and $\ve_0(b)=\ve_1(b)=0$ from the previous 
inequalities. Lemma~\ref{lem:e1 action} then demands that for $\ve_0(b)=\ve_1(b)=0$ to hold, the
diagrams $P$ and $p$ are symmetrical, meaning that the number of columns of the same height
containing only $+$ and only $-$ is  equal, and the number of columns of height $h$ containing 
$\mp$ and  the number of empty columns of height $h-2$ is equal. Therefore $P$ is indeed the 
concatenation of columns as in~\eqref{eq:La pm}. Then it is not hard to see that one can only obtain 
a minimal element from $P$ by the application of a product of appropriate $f(\La_k)$ since otherwise 
there is either a gap in the string amounting to an additional contribution to $\langle c,\ve(b) \rangle$ 
or one of the $-$ in $P$ gives an additional contribution to $\langle c,\ve(b) \rangle$.
The proof that $\vp: B_{\min} \to (P_\cl)_s$ is a bijection is similar.
\end{proof}


\begin{thebibliography}{99}

\bibitem{BFKL:2006}
G.~Benkart, I.~Frenkel, S.-J.~Kang, H.~Lee,
\textit{Level 1 perfect crystals and path realizations of basic representations at $q=0$},
Int. Math. Res. Not.  2006, Art. ID 10312, 28 pp.

\bibitem{C:2001} V.~Chari,
\textit{On the fermionic formula and the Kirillov--Reshetikhin conjecture},
Internat. Math. Res. Notices \textbf{12} (2001) 629--654.

\bibitem{CP:1995}
V.~Chari, A.~Pressley,
\textit{Quantum affine algebras and their representations},
Representations of groups (Banff, AB, 1994),  59--78, CMS Conf. Proc., \textbf{16},
Amer. Math. Soc., Providence, RI, 1995.

\bibitem{CP:1998}
V.~Chari, A.~Pressley,
\textit{Twisted quantum affine algebras},
Comm. Math. Phys. \textbf{196} (1998) 461--476.

\bibitem{FL:2004} 
G.~Fourier, P.~Littelmann,
\textit{Tensor product structure of affine Demazure modules and limit constructions}, 
Nagoya Math. J. \textbf{182} (2006), 171--198. 

\bibitem{FSS:2006}
G.~Fourier, A.~Schilling, M.~Shimozono,
\textit{Demazure structure inside Kirillov-Reshetikhin crystals},
J. Algebra \textbf{309} (2007), 386--404.

\bibitem{Fuchssteiner:1996}
B. Fuchssteiner et al.,
MuPAD User's Manual - MuPAD Version 1.2.2,
John Wiley and sons, Chichester, New York, 1996.
{\tt http://www.mupad.de}

\bibitem{FH:1991}
W.~Fulton, J.~Harris,
Representation Theory. A First Couse,
Graduate Texts in Mathematics, Springer Verlag, 1991.
ISBN: 3-540-97495-4.

\bibitem{H:2007}
D.~Hernandez,
\textit{Kirillov-Reshetikhin conjecture : the general case},
preprint arXiv:0704.2838v1 [math.QA].

\bibitem{HKOTT:2002}
G.~Hatayama, A.~Kuniba, M.~Okado, T.~Takagi, Z.~Tsuboi,
\textit{Paths, crystals and fermionic formulae}, MathPhys odyssey,
2001, 205--272, Prog. Math. Phys., \textbf{23}, Birkh\"auser Boston,
Boston, MA, 2002.

\bibitem{HKOTY:1999}
G.~Hatayama, A.~Kuniba, M.~Okado, T.~Takagi, Y.~Yamada,
\textit{Remarks on fermionic formula},
Recent developments in quantum affine algebras and related topics
(Raleigh, NC, 1998),  243--291, Contemp. Math., \textbf{248},
Amer. Math. Soc., Providence, RI, 1999.

\bibitem{HongKang:2002}
J.~Hong, S.-J.~Kang,
Introduction to Quantum Groups and Crystal Bases,
Graduate Studies in Mathematics, Volume \textbf{42}, American Mathematical Society, 2002.

\bibitem{HN:2006}
D.~Hernandez, H.~Nakajima, 
\textit{Level 0 monomial crystals},
Nagoya Math. J.  \textbf{184}  (2006), 85--153. 

\bibitem{HT:2003}
F.~Hivert, N.~M.~Thi\'ery,
\textit{MuPAD-Combinat, an Open-Source Package for Research in Algebraic Combinatorics},
S\'eminaire Lotharingien de Combinatoire \textbf{51} (2003) [B51z] (70 pp).\newline
{\tt http://mupad-combinat.sourceforge.net/}

\bibitem{JMO:2000}
N.~Jing, K.C.~Misra, M.~Okado,
\textit{$q$-wedge modules for quantized enveloping algebras of classical type},
J. Algebra \textbf{230} (2000), no. 2, 518--539.

\bibitem{Kac:1990}
V.G.~Kac,
Infinite-dimensional Lie algebras. Third edition. Cambridge University Press, Cambridge, 1990. 
xxii+400 pp. ISBN: 0-521-37215-1.

\bibitem{K:1995} 
M.~Kashiwara,
\textit{On crystal bases},
in: Representations of groups (Banff, AB, 1994), 155--197,
CMS Conf. Proc., 16, Amer. Math. Soc., Providence, RI, 1995.

\bibitem{Ka:2002}
M.~Kashiwara,
\textit{On level-zero representation of quantized affine algebras},
Duke Math. J. \textbf{112} (2002), no. 1, 117--195.

\bibitem{KN:1994}
M.~Kashiwara and T.~Nakashima,
\textit{Crystal graphs for representations of the $q$-analogue
of classical Lie algebras},
J. Alg. \textbf{165} (1994) 295--345.

\bibitem{Ko:1999}
Y.~Koga,
\textit{Level one perfect crystals for $B^{(1)}_n,C^ {(1)}_n$, and $D^ {(1)}_n$},
J. Algebra \textbf{217} (1999),  no. 1, 312--334.

\bibitem{KKM:1994} S.-J.~Kang, M.~Kashiwara, K.~C.~Misra,
\textit{Crystal bases of Verma modules for quantum affine Lie algebras},
Compositio Math. \textbf{92} (1994) 299--325.

\bibitem{KKMMNN:1992}
S-J.~Kang, M.~Kashiwara, K.~C.~Misra, T.~Miwa, T.~Nakashima, A.~Nakayashiki,
\textit{Affine crystals and vertex models},
Int.\ J.\ Mod.\ Phys.\ A {\bf 7} (suppl. 1A) (1992) 449--484.

\bibitem{KKMMNN:1992a}
S.-J.~Kang, M.~Kashiwara, K.C.~Misra, T.~Miwa, T.~Nakashima, A.~Nakayashiki,
\textit{Perfect crystals of quantum affine Lie algebras},
Duke Math. J. \textbf{68} (1992), no. 3, 499--607.

\bibitem{KMOY:2006}
M.~Kashiwara, K.~.C.~Misra, M.~Okado, D.~Yamada,
\textit{Perfect Crystals for $U_q(D_4^{(3)})$},
preprint math.QA/0610873.

\bibitem{NS:2005} S.~Naito, D.~Sagaki,
\textit{Construction of perfect crystals conjecturally corresponding to Kirillov--Reshetikhin 
modules over twisted quantum affine algebras},
Comm. Math. Phys.  \textbf{263}  (2006),  no. 3, 749--787.

\bibitem{O:2006}
M.~Okado,
\textit{Existence of Crystal Bases for Kirillov-Reshetikhin Modules of Type $D$},
preprint math.QA/0610874.

\bibitem{OS:2007}
M.~Okado, A.~Schilling,
\textit{Existence of Kirillov-Reshetikhin crystals for nonexceptional types},
preprint arXiv:0706.2224.

\bibitem{S:2002} M. Shimozono,
\textit{Affine type A crystal structure on tensor products of rectangles,
Demazure characters, and nilpotent varieties},
J. Algebraic Combin. \textbf{15} (2002),  no. 2, 151--187.

\bibitem{St:2006}
P.~Sternberg,
\textit{Applications of crystal bases to current problems in representation theory},
PhD thesis, UC Davis 2006 (available at math.QA/0610704).

\bibitem{StS:2006}
P.~Sternberg, A.~Schilling,
\textit{Finite-dimensional crystals $B^{2,s}$ for quantum affine algebras of
type $D_n^{(1)}$},
J. Alg. Combin. \textbf{23} (2006) 317--354.

\bibitem{Y:1998}
S.~Yamane,
\textit{Perfect crystals of $U_q(G^{(1)}_2)$},
J. Algebra  \textbf{210}  (1998),  no. 2, 440--486.

\bibitem{graphviz}
The crystal graph drawings rely on graphviz {\tt http://www.graphviz.org/}

\end{thebibliography}
\end{document}